\documentclass[12pt]{amsart}

\usepackage[all]{xy}
\usepackage{epsfig}
\usepackage{epic}
\usepackage{latexsym}
\usepackage{amsmath}
\usepackage[all]{xy}
\usepackage[mathscr]{euscript}
\usepackage{amssymb}
\usepackage{amscd}
\usepackage{epsfig}

\newtheorem{lemma}{Lemma}[section]
\newtheorem{theorem}[lemma]{Theorem}
\newtheorem{proposition}[lemma]{Proposition}
\newtheorem{corollary}[lemma]{Corollary}

\newcommand{\pf}{\noindent{\em Proof: }}
\newcommand{\epf}{\hfill\hbox{\rule{3pt}{6pt}}\\}

\newcommand{\se}{\!\searrow\!}

\newcommand{\B}{{\mathcal B}}

\newcommand{\V}{{\mathcal V}}
\newcommand{\E}{{\mathcal E}}
\newcommand{\cP}{{\mathcal P}}
\newcommand{\Sg}{{\Sigma}}
\newcommand{\T}{{\mathcal T}}
\newcommand{\M}{{\Sigma^{(\phi)}}}

\newcommand{ \s}{{\Xi}}
\newcommand{\q}{{\cal Q}}

\def\be{\begin{equation}}
\def\q{\end{equation}}

\def\ra{\rightarrow}
\def\lra{\leftrightsquigarrow }
\def\ras{\rightsquigarrow }
\def\Ra{\Rightarrow}

\def\be{\begin{equation}}
\def\ee{\end{equation}}

\def\ra{\rightarrow}

\begin{document}

\title{Blocks and Cut Vertices of the Buneman Graph}

\author{{\bf A.\,W.\,M.\,Dress},
  {\bf K.\,T.\,Huber},
 {\bf J.\,Koolen},
 {\bf V.\,Moulton}}

\address{{\bf Andreas\,W.\,M.\,Dress}\\
Department of Combinatorics and Geometry\\
CAS-MPG Partner Institute and Key Lab for Computational Biology\\
Shanghai Institutes for Biological Sciences\\ Chinese Academy of Sciences,
Shanghai, China\\
and Max Planck Institute for Mathematics in the Sciences\\
D-04103 Leipzig, Germany\\
 email:\,andreas@picb.ac.cn \\
  {\bf Katharina T.\,Huber}\\
 School of Computing Sciences\\
University of East Anglia\\
Norwich, NR4 7TJ, UK\\
 email:\, katharina.Huber@cmp.uea.ac.uk \\
 {\bf Jacobus Koolen}\\
Pohang Mathematics Institute and Department of Mathematics\\
POSTECH, Hyoja-dong, Namgu Pohang\\
790-784 South Korea\\
email: koolen@postech.ac.kr\\
 {\bf Vincent Moulton}\\
School of Computing Sciences\\
University of East Anglia\\
Norwich, NR4 7TJ, UK\\
email:\, vincent.moulton@cmp.uea.ac.uk
}

\date{19.02.2010}
\maketitle

\begin{abstract}
Given a set $\Sg$ of bipartitions of some
finite set $X$ of cardinality at least $2$, one 
can associate to $\Sg$ a canonical 
$X$-labeled graph $\B(\Sg)$, called the Buneman graph.
This graph has several interesting mathematical 
properties --- for example, it is a median network 
and therefore an isometric subgraph of a hypercube. It is commonly used as a tool in studies of DNA 
sequences gathered from populations. 
In this paper, we present some results concerning the
{\em cut vertices} of $\B(\Sg)$, i.e., vertices whose removal 
disconnect the graph, as well as its {\em blocks} or 
$2$-{\em connected components} --- results that yield, in particular, an
intriguing generalization of the well-known 
fact that $\B(\Sg)$ is a tree if and only if any 
two splits in $\Sg$ are compatible.
\end{abstract}

\date{\today}


{\noindent Classification numbers:} 05C05, 05C40, 05C90, 92B10, 92D15

\section{Introduction}

Consider a finite set $X$ of cardinality at least $2$. We 
denote by $\overline A$ the complement $X-A$ for any subset $A$ of $X$. And
we call a bipartition $S=\{A,B\}$
into a proper non-empty subset $A$ of $X$ and its 
complement $B=\overline A$ a {\em split} or, more 
specifically, an $X$-{\em split}.
For any non-empty collection $\Sigma$ of $X$-splits, we define the {\em Buneman graph}
$\B(\Sigma) =\big(V(\Sigma),E(\Sigma)\big)$ 
to be the graph whose vertex set $V(\Sigma)$
consists of all maps $\phi$ from  $\Sigma$ into the power
set $\cP(X)$ of $X$ 
that satisfy, for all $S,S' \in \Sigma$, the
following two conditions

\begin{list}{x}{\leftmargin1cm}
\item[(BG1)] $\phi(S) \in S$, i.e., 
if $S=\{A,B\}$, then $\phi(S)=A$ or $\phi(S)=B$, and

\medskip
\item[(BG2)] $\phi(S) \cap \phi(S') \neq \emptyset$.
\end{list}

\noindent
And we define its edge set $E(\Sigma)$ to consist of all 
those subsets $\{\phi,\psi\}$ of $V(\Sigma)$ for which the 
{\em difference set}
$\Delta(\phi,\psi)$,
defined by 
$$
\Delta(\phi,\psi) : = \{ S \in  \Sigma \,:\, \phi(S)  \neq \psi(S)\}, 
$$ 
has cardinality $1$. 

We also denote by $V^*(\Sigma)$
the superset  of $V(\Sigma)$ 
consisting of {\em all} maps $\phi: \Sigma \to \cP(X)$
that just satisfy (BG1) \big(but not necessarily (BG2)\big). And we define the {\em extended Buneman graph} 
$\B^*(\Sigma)=\big(V^*(\Sigma),E^*(\Sigma)\big)$
to be the (necessarily connected) graph with vertex set
$V^*( \Sigma)$ and edge set $E^*( \Sigma)$ defined
exactly as $E(\Sigma)$ above, yet
with $V( \Sigma)$ replaced by $V^*( \Sigma)$ in its definition. 
Note that $\B^*(\Sigma)$ is clearly isomorphic to an
$|\Sigma|$-dimensional hypercube. 

\begin{figure}[h]
\includegraphics[height=3.5cm, width=7cm]{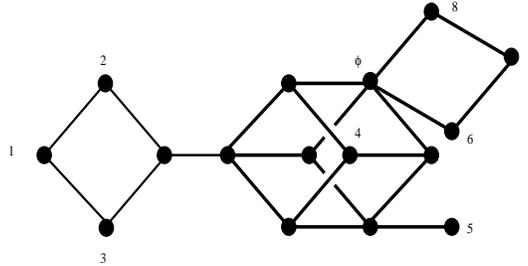}
\caption{For $X :=\{1,2,\dots,8\}$, 
the Buneman graph $\B(\Sigma_8)$ for 
the collection $\Sigma_8$ of $X$-splits given by 
$\Sigma_8 := \{S_{13}, S_{12}, S_{123},
S_{1235}, S_{45}, S_{1234}, S_{67}, S_{78}, S_5\}$
(where,
 for example, $S_{13}$
denotes the $X$-split $\{\{1,3\},\{2,4,5,6,7,8\}\}$).
The vertex $\phi$ is a cut-vertex of $\B(\Sigma_8)$. }
\label{example}
\end{figure}

An example of a Buneman graph is pictured in Figure~\ref{example}. 
Note that some of the vertices of the graph are labelled by elements
in $X$. This arises from a canonical labelling map 
$\varphi_\Sigma:X \ra V(\Sigma): x\mapsto \phi_x$ mapping $X$ 
into the vertex set of $\B(\Sigma)$,
where, for each $x \in X$,  $\phi_x$ denotes the map in $V(\Sigma)$ 
that associates, to any $S \in \Sigma$, 
that subset $S(x)$ of $X$ in $S$ that contains $x$. Clearly, this map is necessarily contained in $V(\Sigma)$. In particular, the cardinality of $V(\Sigma)$ must be at least $2$ for any non-empty collection $\Sigma$ of $X$-splits.

The Buneman graph has appeared in the literature in various guises:
As a co-pair hypergraph in \cite{bar-89,bar-gue-91}, as a 
special type of median graph in e.g.\,\cite{ban-92a} (see 
also \cite{KM99} for a review of median graphs), and in the 
above form in \cite{dre-hub-mou-97}. 
In addition, Buneman graphs and median networks
are regularly used to help analyze viral or mitochondrial 
sequence data gathered from populations 
(see e.g. \cite{BF94,DGL09}).

It is a well-known basic fact essentially established by 
Peter Buneman in 1971 (cf. \cite{B71}) that
$\B(\Sigma)$ is a tree if and only if
any two splits $S,S'$ in $\Sg$ are {\em compatible}, i.e., if and only if, 
for any two splits $S=\{A,B\},S'=\{A',B'\}$ in $\Sg$,
one of the four intersections 
$A\cap A', A\cap B',B\cap A', B\cap B'$ is empty. 

Here, we will generalize this fact in a rather
natural way as follows: Given
{\em any} split system $\Sg$, let $\pi_0(\Sg)$ denote
the set of connected components of the so-called
{\em incompatibility graph} $\Gamma(\Sg)$ of
$\Sg$, i.e., the graph with vertex set
$\Sigma$ whose edge set consists of all
pairs of splits contained in $\Sigma$
that are {\em not} compatible. For example,
for $\Sg_8$ the collection of splits in Figure~\ref{example}, 
$\Gamma(\Sg_8)$ consists of five cliques, {\em viz} $\{S_{13},S_{12}\}$,
$\{S_{123}\}$, $\{S_{1234},S_{1235},S_{45}\}$, 
$\{S_{67},S_{78}\},$
and $\{S_5\}$.

Further,
let $\T(\Sg)$ denote the (by construction
necessarily bipartite) graph with vertex set $\V(\Sg)$
the disjoint union of $V(\Sg)$ and $\pi_0(\Sg)$ and
edge set $\E(\Sg)$ the set of all pairs $\{\Sg_0, \phi\}$
with $\Sg_0\in \pi_0(\Sg), \phi\in V(\Sg)$, and
$\phi(S)\cap A_0\neq \emptyset$ for all $S\in \Sg-\Sg_0$
and all  $A_0\in S_0$ for some --- or, equivalently, every --- split  $S_0\in\Sg_0$. Then, this graph 
$\T(\Sg)=\big(\V(\Sg),\E(\Sg)\big)$ is always a
tree (Theorem~\ref{Xtree}).
For example, the tree $\T(\Sigma_8)$ for the split system $\Sg_8$ referred to in Figure~\ref{example}
is presented in Figure~\ref{construction} (a)
in the last section.

To establish this fact, we will first  
introduce appropriate notations and discuss 
some generalities in Section~\ref{preliminaries}. Then, in 
the next section, we will study the {\em cut vertices} 
of $\B(\Sigma)$, that is, the 
vertices $\phi \in V(\Sigma)$ of $\B(\Sigma)$ for 
which the induced subgraph
$\B^{(\phi)}(\Sigma)$ of $\B(\Sigma)$ with vertex set 
$V(\Sigma)-\{\phi\}$ is disconnected (see e.g. Figure~\ref{example}).  
Remarkably, these vertices can be 
characterized in quite a few equivalent ways (cf. Theorem \ref{cut-vertex}). In  particular, 
given any map $\phi\in V(\Sg)$, one 
can define as many as at least eight significantly 
distinct graphs all of which 
have the property that there is a canonical 
one-to-one correspondence 
between their connected components and the connected 
components of $\B^{(\phi)}(\Sigma)$. So, a map $\phi\in V(\Sg)$ is a cut vertex if and only if any of these graphs is disconnected.

In Section \ref{blocks}, we shall use the cut vertices to provide an explicit description of
the {\em blocks} (or {\em two-connected components}) 
of $\B(\Sigma)$ \cite{D05}, 
that is, the maximal subsets $V$ of $V(\Sigma)$
for which the graph induced on $V$ is connected and does not contain a 
cut vertex. More specifically, we will show that 
there is a canonical one-to-one 
correspondence between $\pi_0(\Sg)$ and the set $Bl(\Sg)$
of blocks of $\B(\Sigma)$ given by associating, to each connected 
component $\Sg_0\in \pi_0(\Sg)$ of $\Gamma(\Sg)$, the set $B(\Sg_0)$ 
consisting of all maps $\phi\in V(\Sg)$ with $\{\Sg_0, \phi\}\in \E(\Sg)$, see 
Theorem~\ref{blockbijection} for details. 
Thus, for example, we see that in Figure~\ref{example}
the Buneman graph has 5 blocks that correspond to the
5 cliques in the associated incompatibility graph. 

And finally, in the last section,
we will establish the above-mentioned generalization of Peter Buneman's 
result from 1971. In addition, we will establish some refinements that will allow us to
associate, to any split system $\Sg$, a ``proper" $X$-tree $T_\Sg$, i.e., 
a triple $(V_\Sg,E_\Sg;\Phi_\Sg)$ consisting 
of a tree with vertex set $V_\Sg $ 
and edge set $E_\Sg $, and a labelling map $\Phi_\Sg: X \ra V_{\Sigma}$ 
such that the degree of any vertex in $V_\Sg$ that is not contained in 
the image of $\Phi_\Sg$ is at least $3$. 
For example, the tree $T_{\Sigma_8}$ for the split 
system $\Sg_8$ in Figure~\ref{example}
is presented in Figure~\ref{construction} (b).

Note that besides providing important structural information 
concerning the Buneman graph (and median networks --- 
see e.g. \cite{BK07}), cut vertices 
have been used to help in the computation of  most 
parsimonious trees for DNA-sequence data 
(see e.g. \cite{BD07,BFA99,S09}). And the blocks of $\B(\Sg)$ determined by them are also closely related to the 
``blobs''  defined in terms of directed graphs described 
by D. Gusfield et al.~in \cite{GB07} (see also e.g. \cite{GB05,HK07}). 

In \cite{DHKM07e}, we will use the results described above to 
show that $\B(\Sg)$ can be described as a push-out
in terms of the data involved in the construction of  $\T(\Sg)$, and 
discuss further potential applications. 

\section{Preliminaries} \label{preliminaries}

In this section, we will review some results 
concerning the Buneman graph that will be needed later
on. Wherever  appropriate, we will refer the reader to 
the relevant literature for proofs of the results that we state. 
We also refer the reader to \cite{D05} for the 
basic terminology of graph theory that we will
use throughout this paper. 

First, we define, for every map $\phi\in V^*(\Sigma)$ and every subset $ \Xi$ of $\Sigma$, 
the map $\phi^ \Xi\in V^*(\Sigma)$ by putting
$$
\phi^\s(S):=\begin{cases}
      & \overline{\phi(S)}\text{ if }S\in \Xi, \\
      & \phi(S)\text{ else, } 
\end{cases}
$$
for every split $S\in\Sigma$. 
Note that, writing $\phi^S$ rather than $\phi^ \Xi$ in case $ \Xi$ 
consists of a single split $S$, only, one has 
$\{\phi,\psi\}\in  E^*(\Sigma)$ for some $\phi,\psi\in V^*(\Sigma)$ 
if and only if
$\psi=\phi^S$ holds for some (necessarily unique) split $S\in \Sigma$, 
{\em viz} the unique split $S=S_{\phi,\psi}$ in $\Delta(\phi,\psi)$. 
Note also that, for a fixed map $\phi \in V^*(\Sigma)$ and 
any other map $\psi\in V^*(\Sigma)$, one clearly  
has $\phi^{\Delta(\phi,\psi)}=\psi$. 
In particular, $\psi$ is completely determined by $\phi$ and 
the difference set $\Delta(\phi,\psi)$. 

Using other notations and arguments, the following result has also (at least implicitly) been shown in, e.g.,  
\cite[Chapter 3.8]{SS03} and \cite[Chapter 5.1]{bar-gue-91}. 
It is also related to \cite[Theorem 1]{GB07}. 

\begin{lemma}
\label{ch3:lemma:buneman:graph:min}
Given a vertex $\phi$ in 
$V(\Sigma)$ and a split $S\in \Sigma$, the following three 
assertions are equivalent:
\begin{enumerate}
\item[(i)] $\phi(S)$ is contained in the set  $\min(\phi[\Sigma])$ 
of inclusion-minimal subsets in 
the image $\phi[\Sigma]$
of $\Sigma$ relative to $\phi$,

\item[(ii)] 
the map $\phi^S\in V^*(\Sigma)$ is, in fact, a vertex in $V(\Sigma)$ 
and thus forms, together with $\phi$, an edge in 
$E(\Sigma)$,

\item[(iii)] there exists some
vertex $\psi$ in 
$V(\Sigma)$
with $\phi(S)\in \min(\phi[\Delta(\phi,\psi)])$.
\end{enumerate}
\end{lemma}
\pf (i) $\Rightarrow$ (ii): Suppose 
$S\in \Sigma$ with $\phi(S) \in \min(\phi[\Sigma])$. 
To see that $\phi^S\in V(\Sigma)$ 
holds, note that, 
by definition, we have 
$\phi^S(S_1)\cap\phi^S(S_2) =\phi(S_1)\cap\phi(S_2) \neq \emptyset$ 
for all $S_1,S_2 \in \Sigma$ that are distinct from $S$. 
Moreover, we have 
$\phi^S(S')\cap\phi^S(S) =\phi(S')\cap\overline{\phi(S)} \neq \emptyset$ 
for all $S' \in \Sigma-\{S\}$ as $\phi(S')\cap \overline{\phi(S)}= \emptyset$ 
would imply $\phi(S')\subset \phi(S)$ in contradiction to 
$\phi(S)\in \min(\phi[\Sigma])$. 

(ii) $\Rightarrow$ (iii): This is trivial: Just put $\psi:=\phi^S$. 
 
(iii) $\Rightarrow$ (i): This follows immediately from the following, 
slightly more general observation:

\begin{lemma}
\label{ch3:lemma:buneman:graph:ideal}Given any two vertices $\phi,\psi$ in 
$\B(\Sigma)$, the set $\phi[\Delta(\phi,\psi)]$ is an {\em ideal} 
 in the partially ordered set $\phi[\Sigma]$, 
that is, 
$A_1,A_2 \in \phi[\Sigma]$, $A_1\in\phi[\Delta(\phi,\psi)]$, and 
$A_2\subseteq A_1$ implies
$A_2\in\phi[\Delta(\phi,\psi)]$. 
\end{lemma}
\pf Indeed, denoting by $S_1$ and $S_2$ the two splits in $\Sigma$ 
with $\phi(S_1)=A_1$ and $\phi(S_2)=A_2$, respectively, we must have
$S_1\in \Delta(\phi,\psi)$. 
Thus, our assumption
$A_2\subseteq A_1=\phi(S_1)$ implies 
$\emptyset=A_2 \cap \overline{\phi(S_1)}=A_2 \cap \psi(S_1)$ and, 
hence, $\psi(S_2)\neq A_2=\phi(S_2)$, i.e., $S_2\in \Delta(\phi,\psi)$, 
as claimed. 

This finishes also the proof 
of Lemma~\ref{ch3:lemma:buneman:graph:min}. \epf

Next, note that, by definition, 
$\B(\Sigma)$ is clearly the induced subgraph of 
$\B^*(\Sigma)$ with vertex set $V(\Sigma)$, and
the graph-theoretical distance $D^*(\psi,\psi')$ 
between any two vertices $\psi$ and $\psi'$
in $\B^*(\Sigma)$ coincides with the cardinality of their difference set $\Delta(\psi,\psi')$. 

It follows that $\B^*(\Sigma)$ is, in particular, a {\em median} graph. I.e., there exists, for any three maps
$\phi_1,\phi_2,\phi_3\in V^*(\Sigma)$, a unique map  $med(\phi_1,\phi_2,\phi_3)$ in $V^*(\Sigma)$, dubbed the {\em median} of $\phi_1,\phi_2,$ and $\phi_3$,
 that 
lies simultaneously on (i) a
shortest path joining $\phi_1$ to $\phi_2$, (ii) a
shortest path joining $\phi_1$ to $\phi_3$, and (iii) a
shortest path joining $\phi_2$ to $\phi_3$. It maps every split $S\in \Sigma$ onto 
\begin{equation}\label{med}
med(\phi_1,\phi_2,\phi_3)(S):=\begin{cases}
       \phi_1(S)&\text{ if }\phi_1(S) \in \{\phi_2(S),\phi_3(S)\}, \\
       \phi_2(S)&\text{ otherwise}. 
\end{cases}
\end{equation}
Parts of the following corollary have also been observed in, e.g., \cite[Section 2.1]{bar-89}, \cite[Proposition 1]{bar-89}, \cite{ban-92a},
and \cite{SS03}:

\begin{corollary}
\label{ch3:lemma:buneman:graph:connected}
The Buneman graph $\B(\Sigma)$ is, for 
every split system $\Sigma$, a connected isometric and median subgraph of the extended
Buneman graph $\B^*(\Sigma)$ associated with $\Sigma$. 
That is:

\noindent
-- The median 
$med(\phi_1,\phi_2,\phi_3)$ of any three maps 
$\phi_1,\phi_2,\phi_3\in V(\Sigma)$ is also contained in 
 $V(\Sigma)$. 
 
 \smallskip
\noindent
 -- The graph-theoretical distance $D(\psi,\psi')$ 
between any two vertices $\psi$ and $\psi'$
in $\B(\Sigma)$ is finite and coincides with the distance $D^*(\psi,\psi')$ between $\psi$ and $\psi'$ in the larger graph $\B^*(\Sigma)$: I.e., 
there exists, for any two vertices $\psi,\psi'
\in V(\Sigma)$, a sequence 
$
\psi_0:=\psi,  \psi_1, \dots, \psi_k :=\psi'
$
of maps in $V(\Sigma)$ of length $k:= |\Delta(\psi,\psi')|$ such that $|\Delta(\psi_{i-1},\psi_i)|=1$ 
holds for all $i=1,\dots,k$. 
 
 More specifically, there exists a canonical one-to-one correspondence 
between the set consisting of all such sequences 
$\psi_0=\psi,  \psi_1, \dots, \psi_{|\Delta(\psi,\psi')|} =\psi'$ and 
the set of all linear orders ``$\preceq$'' defined on 
$\psi[\Delta(\psi,\psi')]$ that extend the partial order of 
$\psi[\Delta(\psi,\psi')]$ defined by set inclusion.

Furthermore, the following holds:
\begin{enumerate}
\item[$(i)$]
Given any three vertices 
$\psi,\psi',\phi\in V(\Sigma)$,
every shortest 
path $\psi_0:=\psi,  \psi_1, \dots, \psi_{|\Delta(\psi,\psi')|} :=\psi'$ 
connecting $\psi$ and $\psi'$ in $\B(\Sigma)$ must pass through $\phi$ if and only if one has
$\Delta(\psi,\psi')=\Delta(\psi,\phi)\cup\Delta(\phi,\psi')$ and 
$\psi(S)\subset \psi(S')$ for every $S\in \Delta(\psi,\phi)$ and 
$S'\in \Delta(\phi,\psi')$. And, conversely,
there exists a shortest path $\psi_0:=\psi,  \psi_1, \dots, 
\psi_{|\Delta(\psi,\psi')|} :=\psi'$ 
connecting $\psi$ and $\psi'$  in $\B(\Sigma)$ not passing through $\phi$ if and only 
if there exist splits $S\in \Delta(\psi,\phi)$ and $S'\in \Delta(\psi',\phi)$
with $\phi(S)\cup \phi(S')\neq X$. 

\item[$(ii)$] 
Any cycle 
$\{\psi_0,\psi_1\},\{\psi_1,\psi_2\},\dots,$ $\{\psi_{n-1},\psi_n:=\psi_0\}$ 
in $E(\Sigma)$ that is not the ``mod-$2$ sum'' -- or symmetric difference -- 
of cycles of smaller length, is of length $4$ 
(as this is easily seen to hold for any median graph). 
\end{enumerate}

\end{corollary}

\pf As, in view of (\ref{med}), there exists, for all $\phi_1,\phi_2,\phi_3\in V^*(\Sigma)$ and all $S,S'\in\Sigma$, an index
$i\in\{1,2,3\}$ with $med(\phi_1,\phi_2,\phi_3)(S)=\phi_i(S)$ and $med(\phi_1,\phi_2,\phi_3)(S')=\phi_i(S')$, it is obvious that $med(\phi_1,\phi_2,\phi_3)\in V(\Sigma)$
holds for any three maps 
$\phi_1,\phi_2,\phi_3\in V(\Sigma)$. 

Further, with $\psi,\psi'\in V(\Sigma)$ as above, select 
a split $S \in \Delta(\psi,\psi')$ 
such that $\psi(S)\in \min\psi[\Delta(\psi,\psi')]$ holds, 
and put $\psi_1:=\psi^S$. 
Clearly, $\psi_1\in V(\Sigma)$,
$\{\psi,\psi_1\}\in E(\Sigma)$, $\Delta(\psi,\psi_1) =\{S\}$, and 
$\Delta(\psi_1,\psi') =\Delta(\psi,\psi')-\{S\}$ and, therefore, 
$\psi_1[\Delta(\psi_1,\psi')] =\psi[\Delta(\psi,\psi')]-\{\psi(S)\}$ also holds. 
Thus, in view of Lemma~\ref{ch3:lemma:buneman:graph:min}, the
first assertion in the statement of the corollary 
follows easily using induction with respect to $|\Delta(\psi,\psi')|$. 

The remaining assertions now follow from this assertion: It implies that, given any three 
vertices $\psi,\psi',\phi\in V(\Sigma)$, the following two claims are equivalent: 

(i) Every shortest path 
$\psi_0:=\psi,  \psi_1, \dots, \psi_{|\Delta(\psi,\psi')|} :=\psi'$ 
connecting $\psi$ and $\psi'$ must pass through $\phi$. 

(ii) 
Both $\Delta(\psi,\psi')=\Delta(\phi,\psi) \cup\Delta(\phi,\psi')$ 
and $\psi(S) \preceq \psi(S')$ hold,
for all $S \in \Delta(\phi,\psi)$ and
$S' \in \Delta(\phi,\psi')$, for  
every linear order ``$\preceq$'' defined on $\psi[\Delta(\psi,\psi')]$ 
that extends the partial order of $\psi[\Delta(\phi,\psi)]$ defined 
by set inclusion. \epf

Next, we associate to every edge 
$e=\{\phi,\psi\}\in E^*(\Sigma)$
its {\em type} $\kappa(e)$ that we define to be the unique split 
$\kappa(e):=S_{\phi,\psi}$ 
in $\Delta(\phi,\psi)$. 
This clearly yields a surjective map $\kappa:E^*(\Sigma)\ra\Sigma$ whose restriction 
to $E(\Sigma)$ we denote by $\kappa_\Sigma$. 

For example, the $4$ horizontal edges of the small cube in the middle of Figure \ref{example} all are of type $S_{1234}$, and the $4$ parallel edges in that cube going from SW to NE are all of type  $S_{1235}$. \\

Clearly,
one has $\Delta(\phi,\psi)=
\{\kappa(\{\phi_{i-1},\phi_i\}): i=1, \dots,|\Delta(\phi,\psi)|\}$ 
for every shortest path 
$\phi_0:=\phi,  \phi_1, \dots, \phi_{|\Delta(\phi,\psi)|} :=\psi$ 
from a vertex $\phi$ to a vertex $\psi$ in $\B(\Sigma)$ or 
$\B^*(\Sigma)$. We claim (see also \cite[Exercise 11, p.64]{SS03})

\begin{lemma} 
The set $\kappa^{-1}_\Sigma(S)$ is, for every split 
$S = \{A,B\} \in \Sigma$, 
a cutset 
of $\B(\Sigma)$ ``inducing'' the split $S$. That is,
removing the edges in $\kappa^{-1}_\Sigma(S)$ from $\B(\Sigma)$ 
yields a subgraph with precisely two connected 
components, one denoted by $\B(\Sigma|A)$ containing all
the vertices $\phi_a$ with $a \in A$, 
and the other one denoted by $\B(\Sigma|B)$ containing all the
vertices $\phi_b$ with $b \in B$. Furthermore, 
the same holds {\rm (mutatis mutandi)} for $\B^*(\Sigma)$ and the subset
$\kappa^{-1}(S)$ 
of $E^*(\Sigma)$.  
\end{lemma}

\pf Let $S = \{A,B\}$ be an arbitrary split in $\Sigma$. It suffices to show that the two subsets
$\B(\Sigma|A) := \{\phi \in V(\Sigma) : \phi(S) = A\}
\supseteq \{\phi_a: a \in A\}$
and $\B(\Sigma|B) := \{\phi \in V(\Sigma) : \phi(S) = B\}\supseteq 
\{\phi_b: b \in B\}$ of $V(\Sigma)$ are (exactly the two) connected
components of the subgraph of $\B(\Sigma)$ obtained by removing the
edges in $\kappa^{-1}_\Sigma(S)$. Indeed, it follows immediately from
the definition of $E(\Sigma)$ that every path in $\B(\Sigma)$ 
from a vertex in $\B(\Sigma|A)$ to a vertex in $\B(\Sigma|B)$ 
must contain at least one edge in $\kappa^{-1}_\Sigma(S)$. 
Moreover, in view of Corollary~\ref{ch3:lemma:buneman:graph:connected},
no shortest path between 
any two vertices in $\B(\Sigma|A)$ \big(and, 
analogously, in $\B(\Sigma|B)$\big),
can pass through an edge in
$\kappa^{-1}(S)$. And, clearly, the same argument 
works just as well for $\B^*(\Sigma)$. 
\epf

Note that restricting the maps $\phi$ 
in $V(\Sigma)$ to a given subset $\Sigma'$ 
of $\Sigma$ clearly induces a surjective 
{\em graph morphism} 
from $\B^*(\Sigma)$  to 
$\B^*(\Sigma')$. That is, it yields a surjective map from $V^*(\Sigma)$ onto 
$V^*(\Sigma')$ that maps every edge $e$ in 
$E^*(\Sigma)$ either onto a single vertex (if $\kappa(e) \in \Sigma-\Sigma'$) in $V^*(\Sigma')$ or onto an 
edge in $E^*(\Sigma')$ (if $\kappa(e) \in \Sigma'$). We now show that
it also induces a surjective graph morphism
${\rm res}_{\Sigma\ra \Sigma'}$ from 
$\B(\Sigma)$  onto $\B(\Sigma')$ (see also \cite{dre-hub-mou-97}):

\begin{lemma}
Given a subset $\Sigma'$ of  a 
split system $\Sigma$ on $X$, and a 
map $\psi \in V(\Sigma')$, there 
exists $($at least$)$ one extension 
$\phi$ of $\psi$ in $V(\Sigma)$, i.e., a 
map $\phi \in V(\Sigma)$ with $\phi\big|_{\Sigma'}=\psi$.

Moreover, the resulting surjective graph morphism 
${\rm res}_{\Sigma\ra \Sigma'}$
contracts every edge $e\in E(\Sigma)$ that 
is of type $\kappa(e)\in \Sigma- \Sigma'$ onto 
a vertex while it maps every edge 
$e\in E(\Sigma)$ of type $\kappa(e)\in \Sigma'$ 
onto an edge $e\in E(\Sigma')$ of the same type.
\end{lemma}

\pf Using induction with respect 
to $|\Sigma|$, we may assume, without loss of generality, that 
$\Sigma=\Sigma'\cup \{S\}$ holds for some single split $S=\{A,B\}$ 
in $\Sigma-\Sigma'$. 
Then, at least one of the two 
extensions $\psi_A,\psi_B$ of $\psi$ in $V^*(\Sigma)$ defined by 
putting $\psi_A(S):=A$ and $\psi_B(S):=B$ must be contained in 
$V(\Sigma)$: Indeed, if $\psi_A\not\in V^*(\Sigma)$ would hold, there would exist some $S' \in \Sigma'$ with 
$\psi_A(S) \cap \psi_A(S')=A\cap \psi(S')= \emptyset$. In this case, however, $\psi_B(S) \cap \psi_B(S'')=B\cap \psi(S'')\neq\emptyset$ 
would hold 
for all  $S''\in \Sigma'$ 
since $B\cap \psi(S'')=\emptyset$ 
for some $S'' \in \Sigma'$ would imply
$\psi(S')\cap \psi(S'')\subseteq  B\cap A=\emptyset$ in 
contradiction to $\psi \in V(\Sigma')$. So, $\psi_B\in V^*(\Sigma)$ must hold in this case. \epf

\begin{corollary}\label{SS'}
Given any two distinct splits $S,S'$ in a split system $\Sg$, 
there exist always two maps $\psi,\psi'\in V(\Sg)$ 
with $S,S'\in\Delta(\psi,\psi')$. 
\end{corollary}
\pf
This follows directly from the last lemma as 
it obviously holds in case $\Sg=\{S,S'\}$.
\epf
\bigskip


\section{Some Graph-Theoretical Observations}
\label{facts}


To establish our main results, we will make use of 
the following simple and purely graph-theoretical observations:

Suppose that $U$ and $V$ are two sets
and that $R\subseteq U\times V$ is a binary relation. 
Let  $\Gamma(R)$ denote the 
bipartite graph with vertex set 
the ``disjoint amalgamation''\footnote{Here, we make use of the 
fact that, according to category theory, the {\em disjoint amalgamation} of any two sets $U$ and $V$ is well defined up to canonical bijection whether or not $U$ and $V$ are disjoint, and can be constructed e.g.~by considering 
the subset of the set $(U \cup V)\times 
\{1,2\}$ consisting of all 
$(w,i) \in (U \times\{1\})\cup (V\times\{2\})$. } 
$U \amalg V$ of $U$ and $V$, and edge set 
$E(R):=\big\{\{u,v\}: u \in U, v\in V, (u,v)\in R\big\}$, 
and define graphs
$\Gamma(R|U):=(U,\big\{\{u_1,u_2\}\in \binom{U}{2} : 
\exists_{v\in V} (u_1,v),(u_2,v)\in R\})$ and 
$\Gamma(R|V):=(V,\big\{\{v_1,v_2\}\in \binom{V}{2} : 
\exists_{u\in U} (u,v_1), (u,v_2)\in R\})$. Let
$\pi_0(R)$, $\pi_0(R|U)$ 
and $\pi_0(R|V)$ denote the connected components of
$\Gamma(R)$, $\Gamma(R|U)$ and $\Gamma(R|V)$, respectively. Then, the following holds:

\begin{lemma}\label{pi0}
Given two sets $U$ and $V$ and a binary relation 
$R\subseteq U\times V$ such that the associated 
bipartite graph $\Gamma(R)$
has no isolated vertices, the embeddings 
$\iota_U:U \ra U \amalg V: u\mapsto u$ and 
$\iota_V:V \ra U \amalg V: v \mapsto v$ induce 
bijections between the  sets $\pi_0(R|U)$ 
and $\pi_0(R|V)$ and the 
set $\pi_0(R)$ of connected components of the graph  
$\Gamma(R)$. I.e., with obvious notational conventions, 
we have a commutative diagram of bijections:

\medskip
\xymatrix{&&
 \pi_0(R|U) 
 \ar@{<->}[rrrr]^-{\pi_0(U\lra V)}
\ar[rrd]_-{\pi_0(\iota_U)}
& && & \pi_0(R|V) \ar@<0.5ex>[lld]^-{\pi_0(\iota_V)}
\\
 &&& &\pi_0(R) 
} 

In particular, given two subsets $A\in \pi_0(R|U)$ 
and $B\in \pi_0(R|V)$ of $U$ and $V$, 
respectively, the following assertions are equivalent:
\begin{itemize}
\item[ (i)] 
$\pi_0(U \lra V)(A) = B$,
\item [(ii)]
$(A\times B) \cap R \neq \emptyset$,
\item [(iii)]
$B=\{v\in V: \exists_{a\in A} (a,v)\in R\}$,
\item [(iv)]
$A=\{u\in U: \exists_{b\in B} (u,b)\in R\}$.
\end{itemize} 
\end{lemma}
\pf All this is quite obvious: The maps  $\pi_0(\iota_U)$ 
and $\pi_0(\iota_V)$ must be surjective as $\Gamma(R)$ 
is supposed to have no isolated vertices. 
And they must be injective because any path 
in $\Gamma(R)$ connecting two vertices 
$u_1,u_2\in U$ (or $v_1,v_2\in V$) gives rise to a path 
connecting these two vertices in $\Gamma(R|U)$ 
(or in  $\Gamma(R|V)$, respectively). \epf  

Now, assume that $U'$ and $V'$ are two further sets, 
and that $\alpha: U'\ra U$ and $\beta: V'\ra V$ are 
two maps such that there exists

\begin{itemize}
\item[(M1)]  for all 
$u_1,u_2\in U$ and $v\in V$ with 
$(u_1,v),(u_2,v) \in R$, some $v'\in V'$ 
with $\big(u_1,\beta(v')\big),\big(u_2,\beta(v')\big) \in R$ and, symmetrically,
\item[(M2)]  
for all 
$u\in U$ and 
$v_1,v_2\in V$ with 
$(u,v_1),(u,v_2) \in R$, some $u'\in U'$ 
with $\big(\alpha(u'),v_1\big),\big(\alpha(u'),v_2\big) \in R$. 
\end{itemize}

Then, defining the binary relations 
$$
R':=R_{\alpha,\beta}:=\{(u',v')\in U'\times V':
\big(\alpha(u'),\beta(v')\big) \in R\},
$$ 
$$
R_\alpha:=\{(u',v)\in U'\times V:
\big(\alpha(u'),v\big) \in R\},
$$ 
and 
$$
R_\beta:=\{(u,v')\in U\times V':
\big(u,\beta(v)\big) \in R\},
$$
it is easily seen that the following holds:
\begin{itemize}
\item[--]
The graph $\Gamma(R_\alpha|V)$ coincides with 
the graph $\Gamma(R|V)$.
\item[--]
The graph $\Gamma(R_\beta|U)$ coincides with 
the graph $\Gamma(R|U)$.
\item[--]
The graph $\Gamma(R_\alpha|U')$ coincides with $\Gamma(R'|U')$
as well as with the graph
induced by $\Gamma(R|U)$ and $\alpha$ on $U'$.
\item[--]
The graph $\Gamma(R_\beta|V')$
coincides with $\Gamma(R'|V')$
as well as with the graph
induced by $\Gamma(R|V)$ and $\beta$ on $V'$. 
\end{itemize}
Furthermore, all the corresponding maps must 
induce bijections on the level of connected components, i.e., we have 

\begin{corollary}\label{pi0cor}
Continuing with the assumptions introduced above
as well as in Lemma~$\ref{pi0}$
and using obvious 
notational conventions, we have 
the following commutative diagram in 
which all maps are bijections:

\medskip
\xymatrix{&&
 \pi_0(R'|U')
 \ar@{<->}[rrrrr]^-{\pi_0(U'\lra V')}
\ar[rrrd]_-{\pi_0(U'\rightsquigarrow R')}
\ar[dddd]_-{\pi_0(U'\rightsquigarrow U)}
& &&& && \pi_0(R'|V') 
\ar[dddd]_-{\pi_0(V'\rightsquigarrow V)}
\ar@<0.5ex>[llld]^-
{\pi_0(V'\rightsquigarrow R')}
\\
 &&& &&
 \pi_0(R') 
\ar@<0.5ex>[ld] _-{\pi_0( R'\rightsquigarrow R_{\alpha}) }
\ar@<-0.5ex>[rd]^-{\pi_0( R'\rightsquigarrow R_{\beta}) }
  \\
 & &&& \pi_0(R_{\alpha}) 
\ar@<0.5ex>[rd] _-{\pi_0( R_{\alpha}\rightsquigarrow R) }
   && \pi_0(R_{\beta}) 
\ar@<0.5ex>[ld] ^-{\pi_0( R_{\beta}\rightsquigarrow R) }
 \\
 &&& &&\pi_0(R) 
 \\
&&
 \pi_0(R|U) 
 \ar@{<->}[rrrrr]^-{\pi_0(U\lra V)}
  \ar[rrru]^-{\pi_0(U\rightsquigarrow R)}
& &&&&& \pi_0(R|V)
  \ar[lllu]_-{\pi_0(V\rightsquigarrow R)}
\\
}
\bigskip
\noindent

In particular, given four subsets $A'\in \pi_0(R'|U')$, 
$A\in \pi_0(R|U)$, $B'\in \pi_0(R'|V')$, 
and $B\in \pi_0(R|V)$ of $U', U, V'$ and $V$, respectively, the following holds:
\begin{itemize}
\item[(i)] 
$\pi_0(U' \lra V')(A') = B'
\iff  (A'\times B')\cap R' \not=\emptyset
\iff B'=\{v'\in V': \exists_{a'\in A'}(a',v')\in R' \}$
\item [(ii)]
$\pi_0(\alpha)(A') = A 
\iff \alpha(A') \cap A\neq \emptyset
\iff\alpha(A') \subseteq A
\iff A=\{u\in U: \exists_{a'\in A',v\in V}(u,v)\in R 
\mbox{ and }(a',v)\in R_{\alpha}\}$
\item [(iii)]
$\big(\pi_0(\beta)\circ\pi_0(U' \lra V')\big)(A')= B 
\iff (\alpha(A')\times B )\cap R \neq \emptyset 
\iff B=\{v\in V:\exists_{a'\in A'}(\alpha(a'),v)\in R \}$
\end{itemize} 
\end{corollary}

\pf  It follows from Lemma~\ref{pi0} that all non-vertical 
maps in the above diagram must be bijections. 
Moreover, if either $\alpha$ or $\beta$ is the identity 
map, at least one of the vertical maps must also be a 
bijection in which case all maps must be bijections. 
So, the general case follows by concatenating the 
diagram obtained for the pair $({\rm id}_U,\beta)$ 
with the diagram obtained for the pair 
$(\alpha,{\rm id}_{V'})$ and applying what we know about 
the individual binary 
relations $R', R_{\alpha},$ and $R$. \epf


\section{Some Characterizations of Cut Vertices}
\label{cutvertices}
%


In this section, we will provide some characterizations
of the cut vertices in the Buneman graph. Some of these closely resemble the characterization of 
cut points in the tight span of a metric space in terms of certain finite graphs given in \cite{DHKM07d}.
More precisely, we will establish the following result:

\begin{theorem}
\label{cut-vertex}
Assume as above 
that $X$ is a finite set, that $\Sigma$ is a 
system of $X$-splits, and that $\phi$ is a map in $V(\Sg)$. Then, the following assertions are equivalent:
\begin{itemize}
\item[(i)] $\phi$ is a cut vertex of $\B(\Sg)$. 

\item[(ii)] There exists a bipartition 
of $\M:=\{ S\in \Sigma:\phi(S)\in \min(\phi[\Sigma])\}$ 
into two disjoint non-empty subsets $\Sg^{(\phi)}_1$ 
and $\Sg^{(\phi)}_2$ such that any two splits $S_1\in\Sg^{(\phi)}_1$ 
and $S_2\in\Sg^{(\phi)}_2$ are compatible. 

\item[(iii)]
$\M$ has a non-empty intersection with at least two distinct connected
components of the incompatibility
graph $\Gamma(\Sg)$ of $\Sg$, as defined in the introduction. 

\item[(iv)]
There exists a bipartition of  $\Sigma$ 
into two disjoint non-empty subsets $\Sg_1$ and $\Sg_2$ such that 
$\phi(S_1)\cup \phi(S_2)=X$ holds for
any two splits $S_1\in\Sg_1$ and $S_2\in\Sg_2$. 

\item[(v)]
There exists a bipartition 
of  $X^{(\phi)}:=\{x\in X:\phi\neq \phi_x\}$ 
into two disjoint non-empty subsets $X_1$ and $X_2$ such that 
$X_1\subseteq\phi(S)$ or $X_2\subseteq\phi(S)$ holds for
all $S\in \Sigma$ or, equivalently, for all $S\in \M$. 

\item[(vi)]
There exists a bipartition of  
$V^{(\phi)}=V^{(\phi)}(\Sigma):=V(\Sigma)-\{\phi\}$ 
into two disjoint non-empty subsets $V_1$ and $V_2$ such that 
$\Delta(\phi,\psi_1)\cap
\Delta(\phi,\psi_2)=\emptyset$ holds for
all $\psi_1\in V_1$ and $\psi_2\in V_2$. 
\end{itemize}
\end{theorem}

\medskip 

To establish Theorem~\ref{cut-vertex}, we apply 
Corollary~\ref{pi0cor} as follows:
With $X, \Sg$, and $\phi$ as in the theorem, 
we put $U:=\Sg$, 
$V:=V^{(\phi)}$, 
$U':=\M$, 
$V':= X^{(\phi)}$.
Further, we denote by $\alpha$ the embedding of $\M$ into $\Sg$ 
and by $\beta$ the map 
$X^{(\phi)}\ra V^{(\phi)}: x\mapsto \phi_x$. And we put 
$R:=R^{(\phi)}:=\{(S,\psi) \in \Sg\times V^{(\phi)}: 
S\in \Delta(\phi,\psi)\}$. 

It is easily checked that
all the requirements needed for applying Corollary~\ref{pi0cor}
are satisfied. For example, the maps $\alpha$ and $\beta$ satisfy
Properties (M1) and (M2): 
Indeed, given two maps $\psi_1,\psi_2\in V^{(\phi)}$, 
and some split $S\in \Sg$ with 
$S\in \Delta(\phi,\psi_1)\cap\Delta(\phi,\psi_2)$, 
there exists a split in $S'\in\M$ with 
$S'\in\Delta(\phi,\psi_1)\cap\Delta(\phi,\psi_2)$ 
as $\psi_1(S)=\psi_2(S)\neq\phi(S) \supseteq \phi(S')$ 
for some $S,S'\in\Sg$ implies 
$\phi(S')\cap  \psi_1(S)= \phi(S')\cap \psi_2(S)=\emptyset$ and, hence,
$ \psi_1(S')=\psi_2(S')\neq \phi(S') $. 
And if $S,S'\in \Delta(\phi,\psi)$ holds for 
some $S,S'\in \Sg$ and $\psi\in V^{(\phi)}$, 
then $ \psi(S)\cap\psi(S')\neq\emptyset$ implies that 
there is some $x\in X^{(\phi)}$ with 
$S,S'\in \Delta(\phi,\phi_x)$ as this 
must hold for any $x\in \psi(S)\cap\psi(S')$.

Now, with $W$ denoting any of 
the sets $\M, X^{(\phi)},\Sg,$ 
and $V^{(\phi)}$ or their cartesian 
products $\M \times X^{(\phi)}, \M \times V^{(\phi)}, \Sg \times X^{(\phi)}$, 
or $\Sg\times V^{(\phi)}$, let $\Gamma_\phi(W)$
denote the corresponding graph with vertex set $W$ whose edge set $E_{\phi}(W)$ is defined in terms of the binary relations 
$R^{(\phi)}, R_\alpha^{(\phi)}, R_\beta^{(\phi)},$ 
and $R^{(\phi)}_{\alpha,\beta}$ -- so, for example, the edge set
$E_{\phi}(\Sigma)$ of $\Gamma_{\phi}(\Sigma)$ is the set 
$$
E_{\phi}(\Sigma)=\{\{S_1,S_2\} \in {\Sigma \choose 2}: \exists_{\psi\in V^{(\phi)}}\,\,
S_1,S_2\in \Delta(\phi,\psi)\}.
$$  
Note that in what comes below, we shall
derive explicit descriptions of 
the edge sets of the graphs 
$\Gamma_{\phi}(\Sg)$, 
$\Gamma_{\phi}(\M)$, 
$\Gamma_{\phi}(X^{(\phi)})$, 
and $\Gamma_{\phi}(V^{(\phi)})$ in (\ref{esg}), 
(\ref{em}), (\ref{exphi}), and (\ref{evphi}), 
respectively. It may be helpful for the reader to take a look at 
these explicit descriptions before proceeding.

In addition,  to further simplify 
notation, put $\pi_\phi(W) := \pi_0\big(\Gamma_{\phi}(W)\big)$
and, for any pair of distinct vertex 
sets $W,W'$ as above, denote by
$\pi_\phi(W' \rightsquigarrow W)$ 
the induced bijection from $\pi_\phi(W')$ onto $\pi_\phi(W)$. 

Then, Corollary~\ref{pi0cor} yields the following  diagram of
canonical bijections:

\medskip
\xymatrix{&&
 \pi_\phi(\M)
 \ar@{<->}[rrrrr]^-{\pi_\phi(\M\lra X^{(\phi)})}
\ar[rrrd]_-{\pi_\phi(\M\rightsquigarrow R^{(\phi)}_{\alpha,\beta})}
\ar[dddd]^-{\pi_\phi(\M \rightsquigarrow \Sg)}
& &&& && \pi_\phi(X^{(\phi)}) 
\ar[dddd]_-{\pi_\phi(X^{(\phi)}\rightsquigarrow V^{(\phi)})}
\ar@<0.5ex>[llld]^-
{\pi_\phi(X^{(\phi)}\rightsquigarrow R^{(\phi)}_{\alpha,\beta})}
\\
 &&& &&
 \pi_\phi(R^{(\phi)}_{\alpha,\beta}) 
\ar@<0.5ex>[ld] 
_-{\pi_\phi( R^{(\phi)}_{\alpha,\beta}\rightsquigarrow R^{(\phi)}_{\alpha}) }
\ar@<-0.5ex>[rd]
^-{\pi_\phi( R^{(\phi)}_{\alpha,\beta}\rightsquigarrow R^{(\phi)}_{\beta}) }
  \\
 & &&& \pi_\phi(R^{(\phi)}_{\alpha}) 
 \ar@<0.5ex>[rd] _-{\pi_\phi( R^{(\phi)}_{\alpha}\rightsquigarrow R^{(\phi)}) }
   && \pi_\phi(R^{(\phi)}_{\beta}) 
    \ar@<0.5ex>[ld] ^-{\pi_\phi( R^{(\phi)}_{\beta}\rightsquigarrow R^{(\phi)}) }
 \\
 &&& &&\pi_\phi(R^{(\phi)}) 
 \\
&&
 \pi_\phi(\Sg) 
  \ar@{<->}[rrrrrr]^-{\pi_\phi(\Sg
  \lra V^{(\phi)})}
  \ar[rrru]^-{\pi_\phi(\Sg\rightsquigarrow R^{(\phi)})}
& &&&&& \pi_\phi(V^{(\phi)})
  \ar[lllu]_-{\pi_\phi(V^{(\phi)}\rightsquigarrow R^{(\phi)})}
\\
} 
\bigskip
\noindent

Now, note that the graph 
$\Gamma_\phi(V^{(\phi)})$ contains the 
induced subgraph $\B^{(\phi)}(\Sigma)$ of
$\B(\Sigma) =
\big(V(\Sigma),E(\Sigma)\big)$ with 
vertex set $V^{(\phi)}$ --- the graph that, by 
definition, is disconnected if and only if $\phi$ 
is a cut vertex (of $\B(\Sigma)$): Indeed, any edge $\{\psi,\psi'\}\in E(\Sigma)$ 
with $\psi,\psi'\neq \phi$ must also 
be an edge in $\Gamma_\phi(V^{(\phi)})$ as $\{\psi,\psi'\}\in E(\Sigma)$ 
implies that either $\Delta(\psi, \phi) \subset\Delta(\psi', \phi)$ 
or $\Delta(\psi', \phi) \subset\Delta(\psi, \phi)$ must hold. 
So, $\psi,\psi'\neq \phi$ implies that
some split $S\in \Delta(\psi, \phi) \cap \Delta(\psi', \phi)$ 
must exist. 

Furthermore, the embedding of $\B^{(\phi)}(\Sigma)$ 
into $\Gamma_\phi(V^{(\phi)})$  induces a bijection 
between the corresponding sets of connected components $\pi_0\big(\B^{(\phi)}(\Sigma)\big)$ and 
$\pi_\phi(V^{(\phi)})$ (the latter being the set 
in the bottom right corner of the above commutative diagram):
Indeed, given any two maps $\psi, \psi'\in V^{(\phi)}$ 
that form an edge in $\Gamma_\phi(V^{(\phi)})$, there 
must exist some split $S \in \Delta(\psi, \phi) \cap \Delta(\psi', \phi) $ 
implying that $S \in 
\Delta(\psi_i, \phi)$ and, hence, that also $\psi_i \neq \phi$ 
must hold for every map $\psi_i$ in any shortest path 
$\psi_0:=\psi, \psi_1, \dots,\psi_k:=\psi'$ from $\psi$ 
to $ \psi'$ in $\B(\Sg)$. Therefore, $\psi$ and $\psi'$ 
must also be contained in the same connected component 
of $ \B^{(\phi)}(\Sigma)$. In consequence, any connected component of $\Gamma_\phi(V^{(\phi)})$ must be contained in and, hence, coincide with  connected component of $ \B^{(\phi)}(\Sigma)$.

Thus, in view of the 
above diagram, 
a map $\phi\in V(\Sg)$ is a cut vertex of $\B(\Sg)$ if and only if 
either one of the eight graphs 
$\Gamma_\phi(W)$ with 
$W=\M, X^{(\phi)},\Sg,V^{(\phi)},  
\M \times X^{(\phi)}, \M \times V^{(\phi)}, 
\Sg \times X^{(\phi)}$, or $\Sg\times V^{(\phi)}$ as above is disconnected.

So, denoting the connected 
component of $\Gamma_\phi(W)$ containing a given 
vertex $w\in W$ by $\Gamma_{\phi}(W,w)$,
Theorem~\ref{cut-vertex} follows immediately from 
the following observations:

\medskip

$\mbox{(i)}\!\!\!\iff\!\! \! \mbox{(iv)}$:
The definition of the edge set $E_\phi(\Sigma)$ implies that 
two distinct splits $S,S'\in \Sigma$ form 
an edge in $E_\phi(\Sigma)$ if and only if there is some 
$x\in \overline{\phi(S)}\cap \overline{\phi(S')}$. 
I.e., we have 
\begin{equation}\label{esg}
E_{\phi}(\Sg) 
= \big\{\{S,S'\}\in \binom{\Sg}{2}:\phi(S)\cup \phi(S') \neq X\big\}.
\end{equation}
So, $\Gamma_\phi(\Sg)$ is disconnected if and only 
if there exists a bipartition of  $\Sigma$ into two 
disjoint non-empty subsets $\Sg_1$ and $\Sg_2$ such that 
$\phi(S_1)\cup \phi(S_2)=X$ holds for
any two splits $S_1\in\Sg_1$ and $S_2\in\Sg_2$. This
establishes the equivalence of (i) and (iv) in Theorem~\ref{cut-vertex}. 

\medskip
It follows also that the incompatibility 
graph $\Gamma(\Sg)$ of $\Sg$ is a subgraph 
of $\Gamma_\phi(\Sigma)$. So, $\Gamma_\phi(\Sigma)$ must be connected for 
every map $\phi\in V(\Sg)$ whenever the 
incompatibility graph $\Gamma(\Sg)$ is connected.
In consequence, $\B^{(\phi)}(\Sigma)$ must 
be $2$-connected in this case, and, more generally,
every connected component in 
$\pi_0(\Sg)=\pi_0\big(\Gamma(\Sg)\big)$ of $\Gamma(\Sg)$ 
must be contained in a connected component 
of $\Gamma_\phi(\Sigma)$. Equivalently, 
every connected component of $\Gamma_\phi(\Sigma)$ 
is a disjoint union of connected components 
of $\Gamma(\Sg)$. The corresponding canonical surjection 
from $\pi_0(\Sg)$ onto 
$\pi_\phi(\Sigma)=\pi_0\big(\Gamma_\phi(\Sigma)\big)$ 
that maps any  connected component $\Gamma(\Sg,S)$ of 
$\Gamma(\Sg)$ containing a split $S\in \Sg$ onto 
the corresponding  connected component $\Gamma_{\phi}(\Sg,S)$ 
of $\Gamma_\phi(\Sg)$ containing $S$
will henceforth be denoted by
$\pi_\phi( \twoheadrightarrow )$.

\medskip
  $\mbox{(i)}\!\!\!\iff\!\!\! \mbox{(ii)}$: 
This follows immediately from the following observation:
\begin{lemma}
\label{incomp}
With $X,\Sg, $ and $\phi$ as above, two splits 
$S,S'\in \Sigma^{(\phi)}$ form an edge in the 
induced subgraph $\Gamma_\phi(\M)$ of $\Gamma_\phi(\Sg)$ 
with vertex set $\M$
if and only if they are incompatible. I.e., we have
\begin{equation}\label{em}
E_{\phi}(\M) 
= \big\{\{S,S'\}\in \binom{\M}{2}:S \mbox{ and } S' 
\mbox{ are incompatible} \big\}.
\end{equation}
\end{lemma}

\noindent {\em Proof of Lemma:} 
Indeed, if two distinct $S,S'\in \Sigma^{(\phi)}$ are 
incompatible, they form an edge in 
$\Gamma_{\phi}(\Sigma)$. 
Conversely, if they 
form an edge in 
$\Gamma_{\phi}(\Sigma)$, i.e., if 
$\overline{\phi(S)}\cap \overline{\phi(S')}\neq \emptyset $ 
holds, they must be incompatible. This follows as  
$\phi(S)\cap \phi(S') \neq \emptyset $ holds 
for any two splits $S,S'\in \Sg$ in view of $\phi\in V(\Sg)$, 
and $\overline{\phi(S)}\cap \phi(S')\neq \emptyset $ holds for any two distinct splits $S\in \M$ and $S'\in \Sg$. \epf
\medskip

$\mbox{(i)}\!\!\!\iff\!\!\! \mbox{(iii)}$: 
Lemma~\ref{incomp}
implies that $\Gamma_\phi(\M)$ can be viewed as 
the induced subgraph of the incompatibility graph $\Gamma(\Sg)$ 
with vertex set $\M\subseteq \Sg$. This implies that 
the embedding $\M\ra \Sg$ induces a well-defined map 
$\pi_\phi(\M) \ra\pi_0(\Sg):
\Gamma_{\phi}(\M,S) \mapsto \Gamma(\Sg,S)$ 
that we denote, for short, by
$\pi_\phi(\rightarrowtail)$.
Moreover, the composition 
$$
\pi_\phi( \twoheadrightarrow )  \circ \pi_\phi(\rightarrowtail): 
\pi_\phi(\M)\ra \pi_0(\Sg)\ra \pi_\phi(\Sg)
$$ 
maps, for any split $S\in \M$, the connected 
component $\Gamma_\phi(\M,S)$
onto the connected component $\Gamma_\phi(\Sg,S)$ and, hence, 
coincides with the bijection 
$\pi_\phi(\M \rightsquigarrow \Sg):
\Gamma_{\phi}(\M,S) \mapsto \Gamma_{\phi}(\Sg,S)$ 
induced by the embedding $\alpha$ of $\M$ into $\Sg$. 
In turn, this implies 

\begin{lemma}\label{M}
The map $\pi_\phi(\rightarrowtail)$ is always injective, 
that is, any two splits $S,S'\in \M$ can be connected by a 
sequence of pairwise incompatible splits in $\M$ whenever 
they can be connected by such a sequence of pairwise 
incompatible splits in $\Sg$.
\end{lemma}

This clearly establishes that
$\mbox{(i)}\!\!\!\iff\!\!\! \mbox{(iii)}$ holds, as claimed.

\medskip
$ \mbox{(i)}\!\!\!\iff\!\!\! \mbox{(v)}\!\!\!\iff\!\!\! \mbox{(vi)}$:  
We now establish the equivalence of the first 
and the last two assertions in Theorem~\ref{cut-vertex}: It suffices to note that two elements $x_1,x_2$ in $X^{(\phi)}$ 
are {\em not} connected by an edge in $\Gamma_\phi(X^{(\phi)})$ 
if and only if $x_1\in\phi(S)$ or $x_2\in\phi(S)$ holds for
all $S\in \Sigma$ or, equivalently, for all $S\in \M$. I.e., we have
\begin{eqnarray}\label{exphi}
&& E_{\phi}(X^{(\phi)}) \\
&&= \big\{\{x_1,x_2\}\in \binom{X^{(\phi)}}{2}:
\exists_{S \in \Sg \mbox{ \tiny \em (or $\M$)}} \,\,\,
x_1,x_2 \not\in \phi(S)\big\}\nonumber \\ 
&&= \big\{\{x_1,x_2\}\in \binom{X^{(\phi)}}{2}: 
D(\phi_{x_1},\phi_{x_2})<D(\phi_{x_1},\phi)+D(\phi,\phi_{x_2})  
\big\}.\nonumber \end{eqnarray}
Moreover, two maps $\psi_1,\psi_2$ in $V^{(\phi)}$ 
are {\em not} connected by an edge in 
$\Gamma_\phi(V^{(\phi)})$ if and only 
if $\phi(S)=\psi_1(S)$ or $\phi(S)=\psi_2(S)$ holds for
all $S\in \Sigma$ or, equivalently, for all $S\in \M$. I.e., we have
\begin{eqnarray}\label{evphi}
&&E_{\phi}(V^{(\phi)})\\
&&  = \big\{\{\psi_1,\psi_2\}\in \binom{V^{(\phi)}}{2}:
D(\psi_1,\psi_2)<D(\psi_1,\phi)+D(\phi,\psi_2)  \big\}\nonumber\\
&&= \big\{\{\psi_1,\psi_2\}\in \binom{V^{(\phi)}}{2}:
\exists_{S \in \Sg \mbox{ \tiny (or $\M$)}}\,\,\,
\psi_1(S)=\psi_2(S) \neq \phi(S)\big\}.\nonumber
\end{eqnarray}
This finishes the proof of Theorem \ref{cut-vertex}. \epf

\medskip
It is worth noting in this context that 
our approach implies also that, given a map $\phi\in V(\Sg)$, the following assertions are equivalent:
\begin{itemize}
\item[(i)] $\phi$ is a cut vertex.
\item[(ii)] 
There exists a bipartition 
$\{\Sg_1,\Sg_2\}$ 
of $\Sigma$ and a bipartition 
$\{V_1,V_2\}$ 
of $V^{(\phi)}$ such 
that $\phi(S)=\psi(S)$ holds for all $S\in \Sg_1$ 
and $\psi\in V_2$, and for all $S\in \Sg_2$ and 
$\psi\in V_1$.
\item[(iii)]
There exists a bipartition 
$\{\Sg_1,\Sg_2\}$ 
of $\M$ and a bipartition 
$\{X_1,X_2\}$ of $X^{(\phi)}$ such 
that $x\in \phi(S)$ holds for all $S\in \Sg_1$ 
and $x \in X_2$, and for all $S\in \Sg_2$ and $x\in X_1$.
\end{itemize}
For example, for the map $\phi$ 
in Figure~\ref{example}, the bipartition of $\Sigma_8$ is
$\big\{\{S_{67},S_{78}\},\Sigma - \{S_{67},
S_{78}\}\big\}$. 
Moreover, we have
$\Sigma^{(\phi)} = \{S_{1235},S_{1234}, S_{45},S_{78},
S_{67}\}$,
and the corresponding bipartition is given by 
$\big\{\{S_{67},S_{78}\},\Sigma^{(\phi)} - 
\{S_{67},S_{78}\}\big\}$.

It is also worth noting that the various images of the
connected components in the sets $\pi_\phi(...)$ 
relative to the respective bijections considered 
above can be described as follows. 
For $\phi,\psi \in V^*(\Sigma)$, we put 
$$
\Delta_{\min}(\psi|\phi):=\{S\in \Delta(\phi,\psi): 
\psi(S)\in \min(\psi[\Delta(\phi,\psi)])\},
$$ 
where $\min(\psi[\Delta(\phi,\psi)])$ denotes
the set of (inclusion-)minimal subsets in 
the image $\psi[\Delta(\phi,\psi)]$
of $\Delta(\phi,\psi)$ relative to $\psi$. 

\begin{proposition}\label{1-1}
Given any four connected 
components $\Sigma_0' \in \pi_{\phi}(\Sigma^{(\phi)})$,
$\Sigma_0 \in \pi_{\phi}(\Sigma)$,
$X_0 \in \pi_{\phi}(X^{(\phi)})$, and 
$V_0 \in \pi_{\phi}(V^{(\phi)})$, the following holds:

\begin{itemize}
\item[(i)]  
$\pi_\phi(\M\ras\Sg)(\Sg_0')=\Sg_0 
\iff \Sg_0'\subseteq \Sg_0
\iff \Sg_0'= \Sg_0\cap \M \iff \Sg_0=
\{S\in\Sg:\phi(S)\cup\phi(S')\neq X\mbox{ for some } S'\in \Sg_0'\}$. 

\medskip

\item[(ii)] 
$\pi_\phi(\M\leftrightsquigarrow X^{(\phi)})(\Sg_0')=X_0
\iff X_0=\{x\in X^{(\phi)}:
\Delta_{\min}(\phi|\phi_x)\subseteq
\Sg_0' \} 
\iff 
X_0=\{x\in X^{(\phi)}:
\Delta_{\min}(\phi|\phi_x)\cap
\Sg_0' \neq \emptyset\}
\iff \Sg_0'=\bigcup_{x \in X_0} \Delta_{\min}(\phi|\phi_x)$. 

\medskip

\item[(iii)] 
$\pi_\phi(\M\ras V^{(\phi)})(\Sg_0')=V_0
\iff V_0=\{\psi\in V^{(\phi)}:
\Delta_{\min}(\phi|\psi)\subseteq
\Sg_0' \} \iff V_0=\{\psi\in V^{(\phi)}:
\Delta_{\min}(\phi|\psi)\cap
\Sg_0' \neq \emptyset\}
\iff V_0\supseteq\{\phi^S: S\in\Sg_0' \}
\iff \Sg_0'=\bigcup_{\psi\in V_0} \Delta_{\min}(\phi|\psi) \iff \Sg_0'=
\{S\in\M: \phi^S \in V_0\}$. 

\medskip

\item[(iv)] 
$\pi_\phi(\Sg\ras X^{(\phi)})(\Sg_0)=X_0
\iff X_0=\{x\in X^{(\phi)}:
\Delta(\phi,\phi_x)\subseteq
\Sg_0 \}\iff X_0=
\big\{x\in X^{(\phi)} :\Sigma_0\cap \Delta(\phi,\phi_x)
\not=\emptyset\big\}
 \iff \Sg_0=\bigcup_{x \in X_0} \Delta(\phi,\phi_x)$. 

\medskip

\item[(v)] 
$\pi_\phi(\Sg\leftrightsquigarrow V^{(\phi)})(\Sg_0)=V_0
\iff V_0=\{\psi \in V^{(\phi)}:
\Delta(\phi,\psi)\subseteq
\Sg_0 \} 
\iff V_0=
\big\{\psi\in V^{(\phi)}: \Delta(\phi,\psi)\cap \Sigma_0\not=\emptyset\}
\iff \Sg_0=\bigcup_{\psi \in V_0} \Delta(\phi,\psi)$. 

\medskip

\item[(vi)] 
$\pi_\phi(X^{(\phi)}\ras V^{(\phi)})(X_0)=V_0
\iff V_0=\{\psi \in V^{(\phi)}:
D(\psi,\phi_x)<D(\psi,\phi) + D(\phi,\phi_x) \mbox{ for some } x\in X_0 \}
\iff V_0=\{\psi \in V^{(\phi)}:
 \Delta(\phi,\phi_x)\cap \Delta(\phi,\psi)
\not=\emptyset\mbox{ for some } x\in X_0 \}
 \iff X_0=\{x\in X:\phi_x\in V_0\}$. 
\end{itemize}
\end{proposition}

These assertions follow quite easily from our definitions. 
We leave their simple and straight-forward 
(yet sometimes a bit laborious) verification to the 
interested reader. 

Note also that, continuing with the assumptions 
and notations of Theorem~\ref{cut-vertex}, our 
analysis implies the following corollary:

\begin{corollary}\label{2-connected}
The graph $\B(\Sg)$ is $2$-connected if and only if the incompatibility
graph $\Gamma(\Sg)$ is connected. 

More generally, given any two maps
$\psi,\psi'\in V(\Sg)$, one has 
$\Gamma_{\phi}(V^{(\phi)},\psi)$\,\,\,\,
$= \Gamma_{\phi}(V^{(\phi)},\psi')$ 
for all $\phi\in V(\Sg)-\{\psi,\psi'\}$ if and 
only if $\Delta(\psi,\psi')\subseteq \Sg_0$ holds 
for some connected component $\Sg_0$ of the 
incompatibility graph $\Gamma(\Sg)$
of $\Sg$. So, conversely, 
two maps $\psi,\psi'\in V(\Sg)$ are 
separated by some cut vertex $\phi\in V(\Sg)$ 
if and only if $\Delta(\psi,\psi')$ has a 
non-empty intersection with at least two distinct 
connected components of $\Gamma(\Sg)$. 
\end{corollary}

\pf We have already seen above that 
$\B(\Sg)$ cannot contain a cut vertex 
in case $\Gamma(\Sg)$ is connected.

Conversely, if $\Gamma(\Sg)$ is not connected, 
one may choose any two splits $S,S'$ in distinct 
connected components of $\Gamma(\Sg)$ and then, 
according to Corollary~\ref{SS'}, 
two maps $\psi,\psi'\in V(\Sg)$ with $S,S'\in\Delta(\psi,\psi')$. Then, 
given any shortest path 
$\psi_0:=\psi, \psi_1, \dots,\psi_k:=\psi'$ from 
$\psi$ to $ \psi'$ in $\B(\Sg)$, 
there must exist some $i$ in $\{1,\dots,k-1\}$ 
such that the two splits in the one-split sets 
$\Delta(\psi_i,\psi_{i-1})$ and $\Delta(\psi_i,\psi_{i+1})$ 
are in distinct connected components of $\Gamma(\Sg)$. Thus,  
$\psi_i$ must be a cut vertex in $\B(\Sg)$. 

Moreover, 
if $\Delta(\psi,\psi')\subseteq \Sg_0$ holds 
for some connected component $\Sg_0$ of the 
incompatibility graph $\Gamma(\Sg)$
and $\phi$ is any map in $V(\Sg)-\{\psi,\psi'\}$, 
we must have $\Gamma_{\phi}(V^{(\phi)},\psi)=\Gamma_{\phi}(V^{(\phi)},\psi')$:  
Indeed, Proposition \ref{1-1} (v) implies that $\Gamma_{\phi}(V^{(\phi)},\psi)\neq \Gamma_{\phi}(V^{(\phi)},\psi')$ holds for some maps $\psi,\psi' \in V^{(\phi)}$ if and only if 
the two sets $\Delta(\psi,\phi)$ and
$\Delta(\phi,\psi')$ are contained in two distinct connected components of $\Gamma_\phi(\Sg)$. 
In turn, this implies that 
$\Delta(\psi,\psi')$ must coincide
with the disjoint union of $\Delta(\psi,\phi)$ and
$\Delta(\phi,\psi')$. So, $\Delta(\psi,\psi')$ cannot be contained in a single connected component of $\Gamma_\phi(\Sg)$ and, hence, even less in a single connected component of $\Gamma(\Sg)$ in this case. Thus, given any two maps
$\psi,\psi'\in V(\Sg)$, 
$\Gamma_{\phi}(V^{(\phi)},\psi)=\Gamma_{\phi}(V^{(\phi)},\psi')$ must hold indeed for any map
 $\phi$ in $V(\Sg)-\{\psi,\psi'\}$ in case 
$\Delta(\psi,\psi')\subseteq \Sg_0$ holds 
for some connected component $\Sg_0$ of the 
incompatibility graph $\Gamma(\Sg)$.

Conversely, if 
$\Gamma_\phi(V^{(\phi)},\psi)=\Gamma_{\phi}(V^{(\phi)},\psi')$ 
holds for all $\phi\in V(\Sg)-\{\psi,\psi'\}$, 
we must have  $\Delta(\psi,\psi')\subseteq \Sg_0$ for 
some connected component $\Sg_0$ of the incompatibility 
graph $\Gamma(\Sg)$. Indeed, choosing as above any path 
$\psi_0:=\psi, \psi_1, \dots,\psi_k:=\psi'$ from $\psi$ 
to $ \psi'$ in $\B(\Sg)$, 
there would otherwise exist some $i$ in $\{1,\dots,k-1\}$ 
such that the two splits in the one-split sets 
$\Delta(\psi_i,\psi_{i-1})$ and $\Delta(\psi_i,\psi_{i+1})$ 
are in distinct connected components of $\Gamma(\Sg)$.
In turn, this implies that, for $\phi:=\psi_i$, we would have  
$\Gamma_\phi(V^{(\phi)},\psi_{i-1}) \neq
\Gamma_\phi(V^{(\phi)},\psi_{i+1})$ 
in view of Proposition~\ref{1-1}(iii), and 
the fact (cf. Lemma~\ref{M}) that distinct 
connected components of $\Gamma(\Sg)$ intersect $\M$ in distinct connected components of $\Gamma_\phi(\M)$ and, hence, $\Gamma_\phi(V^{(\phi)},\psi)=\Gamma_\phi(V^{(\phi)},\psi_{i-1}) \neq
\Gamma_\phi(V^{(\phi)},\psi_{i+1})=\Gamma_{\phi}(V^{(\phi)},\psi')$, a contradiction. 
\epf

To conclude this section, we note that, essentially by definition, $(S,\phi^S)\in R^{(\phi)}$ holds for every $\phi\in V(\Sg)$ and all $S\in \M$. Thus, the obviously well-defined and injective 
map $\gamma: \M\ra V^{(\phi)}: S\mapsto \phi^S$ 
necessarily induces a map $\pi_0(\gamma)$ from $\pi_{\phi}(\Sigma^{(\phi)})$ 
into $\pi_{\phi}(V^{(\phi)})$ that must
coincide with the bijection $\pi_\phi(\M \ras V^{(\phi)})$. So, we have a diagram of bijections all of which are ``induced''
by naturally defined maps between 
the corresponding vertex sets:

\medskip
\xymatrix{&&
 \pi_\phi(\M) 
 \ar@{->}[rrrrdd]^-{\pi_0(\gamma)}
\ar[dd]_-{\pi_\phi(\M\ras\Sg)}
& && & \pi_\phi(X^{(\phi)}) 
\ar[dd]^-{\pi_\phi(X^{(\phi)}\ras V^{(\phi)})}
\\
 &&& &
\\
&&
 \pi_\phi(\Sg) 
& && & \pi_\phi(V^{(\phi)}) 
\\
} 
\bigskip

\noindent
So, once established, these bijections could also have been used to 
define all the other bijections between the sets 
$\pi_\phi(\M) , \pi_\phi(\Sg) ,\pi_\phi(X^{(\phi)})$, and 
$\pi_\phi(V^{(\phi)})$ which, however, we feel would have 
led to a much less transparent and natural approach.

\section{Blocks of the Buneman graph}\label{blocks}

For a collection $\Sigma$ of $X$-splits, we denote
by $cut(\Sigma)$ the set of all cut vertices of
$\B(\Sigma)$ and, as mentioned in the introduction, by $Bl(\Sigma)$
the set of all of its blocks.
In this section, we will describe a canonical bijection between 
$Bl(\Sigma)$ and $\pi_0(\Sigma)$, the set  
of connected components of the incompatibility graph of $\Sigma$.

To this end, note first that, given any two distinct compatible $X$-splits $S$
and $S'$, there exists a unique subset $A\in S$ that we
denote by $A(S\se S')$ such that $A\cap A'\neq \emptyset$ 
and $A\cap B'\neq \emptyset$ or, equivalently, 
$A'\subset A$ or $B'\subset A$ holds for the two subsets $A',B'$ in $S'$. 
Clearly,  given $A\in S$ and $A'\in S'$, we have 

\begin{itemize}
\item[(4.1)]
$A\cup A' = X \iff  A= A(S\se S')$ and $A'= A(S'\se S)$.
\medskip

\item[(4.2)]
$A'\subset A\iff A= A(S\se S')$ and $A'\neq A(S'\se S)$.
\medskip

\item[(4.3)]
$A\subset A'\iff A\neq A(S\se S')$ and $A' = A(S'\se S)$.
\medskip

\item[(4.4)]
$A'\cap A=\emptyset\iff A\neq A(S\se S')$ and $A'\neq A(S'\se S)$.
\medskip

\item[(4.5)]
If $\phi\in V(\Sg)$
and $S\in \M$, then
$\phi(S')=A(S'\se S)$ must hold for every 
split $S'\in \Sg$ that is compatible with  $S$.
\medskip

\pf Just apply (4.2) and (4.4) 
to $A:=\phi(S)$ and $A':=\phi(S')$, noting that neither 
$A'\subset A$ nor $A'\cap A=\emptyset$ can hold.\epf 
\end{itemize}

Similarly, suppose that 
$S,S',S''$ are three distinct $X$-splits 
such that $S$ is compatible to 
$S'$ and to $S''$. Then,
as both $A(S'\se S)$ and $A(S''\se S)$ properly 
contain one of the two sets in $S$, we have
\medskip

\begin{itemize}
\item[(4.6)]
$A(S'\se S)\cap  A(S''\se S) \neq \emptyset$, and
\medskip

\item[(4.7)] 
$A(S\se S')=A(S\se S'')$ whenever $S'$ and $S''$ are
incompatible.
\medskip

\pf $A(S\se S')\neq A(S\se S'')$ 
would imply $A(S\se S'')  \subset A(S'\se S)$ 
\big(apply (4.3) with $A :=A(S\se S'')$\big) and hence, by symmetry, also 
$A(S\se S') \subset A(S''\se S)$.
Thus $X =A(S\se S'')\cup A(S\se S') \subseteq A(S'\se S)\cup
A(S''\se S)$, implying that $S'$ and $S''$ 
must be compatible, a contradiction.\epf
\end{itemize}

Thus, given any connected component $\Sg_0$ of the 
incompatibility graph $\Gamma(\Sg)$ of $\Sg$ and any 
split $S\in \Sg -\Sg_0$, the two subsets $A(S\se S')$ and $A(S\se S'')$ 
in $S$ are well-defined and coincide for any 
two splits $S',S''\in \Sg_0$. In consequence, we will also 
write $A(S\se \Sg_0)$
for this subset in $S$.

Now, with $\Sigma$ a collection of $X$-splits 
as above, consider a connected component 
$\Sigma_0 \in \pi_0(\Sigma)$ of the 
incompatibility graph $\Gamma(\Sg)$ of $\Sg$ 
and associate, to any map $\phi \in V^*(\Sigma_0)$, 
the map $\widetilde{\phi}: \Sigma \to \cP(X)$ defined by putting
$$
\widetilde{\phi}(S) := 
\left\{\begin{array}{cc}
\phi(S) & \mbox{ if } S\in \Sigma_0,\\
A(S\se \Sigma_0)& \mbox{ else.}
\end{array}
\right.
$$ 
We claim

\begin{theorem}\label{blockbijection}
With $\Sigma$ and $\Sigma_0$ as above, 
we have $\widetilde{\phi}\in V(\Sigma)$ for every map $\phi \in V(\Sigma_0)$.  
Moreover,
the corresponding embedding
$$
\Theta=\Theta_{\Sg_0}:
V(\Sigma_0)\ra V(\Sigma):
\phi\mapsto  \widetilde{\phi}
$$
is actually an isometry 
from
$V(\Sigma_0)$ into $V(\Sigma)$ that maps $V(\Sigma_0)$
onto a block 
$$
B(\Sg_0)=B_\Sg(\Sg_0):= \{ \widetilde{\phi} \,:\, \phi \in V(\Sigma_0)\}
$$ 
of $\B(\Sigma)$. 
Furthermore, we have
\begin{eqnarray}\label{B0}
B(\Sigma_0) 
& = & \{\phi\in V(\Sigma): \forall_{S\in \Sigma-\Sigma_0} 
\phi(S)=A(S\se \Sigma_0)\}\\
& = & \{\phi\in V(\Sigma): \M\cap\Sigma_0 \neq \emptyset\}\nonumber
\end{eqnarray}
and, hence, also 
\begin{eqnarray}\label{B00}
B(\Sigma_0) 
 =\{\phi\in V(\Sigma): \Delta(\phi,\phi_0)\subseteq \Sigma_0)\}
 \end{eqnarray}
 for any map $\phi_0\in B(\Sigma_0)$.
And there exist a $($necessarily unique$)$ ``gate'' $\phi_{\Sigma_0}\in B(\Sigma_0)$ for any map $\phi\in V(\Sg)$, i.e., 
a map in $B(\Sigma_0)$ 
with $D(\phi, \psi)=D(\phi, \phi_{\Sigma_0}) + D( \phi_{\Sigma_0}, \psi)$ 
for all $\psi\in B(\Sigma_0)$, which is given by
$$
 \phi_{\Sigma_0}(S):=\begin{cases}
      \,\phi(S) \text{ if } S\in \Sg_0, \\
      \,A(S\se \Sigma_0) \text{ otherwise}
\end{cases}
$$
and must necessarily be contained in $cut(\Sg)$ in case $\phi \not\in B(\Sigma_0)$ holds.

In particular, given any two distinct connected components $\Sigma_0,\Sigma_1\in \pi_0(\Sg)$ of the incompatibility graph of $\Sg$, the map defined by 
$$
\phi_{\Sg_0|\Sg_1}: \Sigma\ra \cP(X):S \mapsto \begin{cases}
      \, A(S\se \Sigma_1)\text{  if  } S\in \Sigma_0, \\
       \,A(S\se \Sigma_0)\text{ otherwise,}
\end{cases}
$$
is the unique gate $\phi_{\Sg_0}\in B(\Sigma_0)$ of all maps $\phi \in B(\Sigma_1)$. And any map in the intersection $B(\Sigma_0)\cap B(\Sigma_1)$ of the blocks $B(\Sigma_0)$ and $B(\Sigma_1)$, if there is any such map, must be a cut vertex of $\B(\Sg)$. 
Furthermore, the following assertions all are equivalent: 
\begin{itemize}
\item[($B\cap B:1$)] There exists some 
$\phi \in cut(\Sigma)$ with $\phi\in B(\Sigma_0)\cap B(\Sigma_1)$.
\smallskip
\item[$(B\cap B:2)$]  $B(\Sigma_0)\cap B(\Sigma_1)$is non-empty.
\smallskip
\item[$(B\cap B:3)$]  $B(\Sigma_0)\cap B(\Sigma_1)$ is a one-vertex set.
\smallskip
\item[$(B\cap B:4)$]  $A(S\se \Sigma_0)=A(S\se \Sigma_1)$ 
holds for all $S\in \Sg-(\Sigma_0\cup\Sigma_1)$. I.e.,  there exists no split $S=\{A,B\}\in \Sg-(\Sigma_0\cup\Sigma_1)$ with 
$A\cap A_0=\emptyset$ for some set $A_0$ in 
some split $S_0\in \Sigma_0$ and $B\cap B_1=\emptyset$ 
for some set $B_1$ in some split $S_1\in \Sigma_1$.
\end{itemize} 
And in case all of these assertions hold, the 
map $\phi_{\Sg_0|\Sg_1}\in B(\Sg_0)$ coincides with the correspondingly defined map $\phi_{\Sg_1|\Sg_0}$, and it is the unique map in $B(\Sigma_0)\cap B(\Sigma_1)$. 

Finally, mapping each connected component $\Sg_0$ of the
incompatibility graph of $\Sg$ onto
the associated block $B(\Sg_0)$ induces a canonical bijection
$$
\Psi=\Psi_{\Sg}: \pi_0(\Sigma) \to Bl(\Sigma): \Sigma_0 \mapsto B(\Sg_0)
$$
from the set $\pi_0(\Sigma)$ of connected 
components of $\Gamma(\Sg)$ onto the set $Bl(\Sigma)$
of all blocks of $\B(\Sigma)$.
\end{theorem}

\pf
To simplify the exposition of the proof, we will
present it as a series of 12 observations:
(i) It follows immediately from the definitions and
Assertion~(4.6) that $B(\Sigma_0)$ is a subset of 
$V(\Sigma)$. It is also obvious that $\Theta$ is 
an isometry as, by definition, the even stronger assertion
$$
\Delta(\widetilde{\phi}, \widetilde{\psi})=\Delta(\phi, \psi)
$$ 
apparently holds for all $\phi,\psi\in V(\Sg_0)$. Hence,  the subgraph induced by $\B(\Sigma)$ on
$B(\Sigma_0)$ is an isometric subgraph of $\B(\Sigma)$.
\medskip

(ii) We clearly have 
$$
B(\Sigma_0)=\{\phi\in V(\Sigma): \forall_{S\in \Sigma-\Sigma_0}
\phi(S)=A(S\se \Sigma_0)\}
$$ 
as $\phi(S)=A(S\se \Sigma_0)$ holds, by definition, 
for all $\phi\in B(\Sigma_0)$ and all
$S\in \Sigma-\Sigma_0$. Conversely, given any 
map $\phi\in V(\Sigma)$ with $\phi(S)=A(S\se \Sigma_0)$ 
for all $S\in \Sigma-\Sigma_0$, we have 
$\psi:=\phi|_{\Sigma_0}\in V(\Sigma_0)$ and 
$\widetilde{\psi}=\phi$. 
\medskip

(iii) Further, if $\phi(S')=A(S'\se \Sigma_0)$ holds
for all $S'\in \Sigma-\Sigma_0$, $\phi(S') \cap A \neq \emptyset$ must hold 
for all $S'\in \Sigma-\Sigma_0$ and subsets $A$ 
in any split $S\in \Sigma_0$. Thus, if $\phi=\widetilde{\psi}$ holds for some map $\psi\in V(\Sigma_0)$, $S\in \M$ must hold for any $S\in \Sigma_0$ 
for which $\phi(S)=\psi(S)$ is a minimal subset 
in $\psi[\Sg_0]=\{\psi(S'): S'\in \Sg_0\}$. 
So, we have indeed
$
B(\Sigma_0)\subseteq\{\phi\in V(\Sigma): \M\cap\Sigma_0 \neq \emptyset\},
$
as required.

\medskip

(iv) Conversely, given any 
$\phi\in V(\Sigma)$ with $\M\cap\Sigma_0 \neq \emptyset$, 
Assertion~(4.5) implies that one must have $\phi(S')=A(S'\se S)$
for any split $S\in \M$ and any split 
$S'\in \Sg$ that is compatible with $S$.
Hence, $\phi(S')=A(S'\se \Sg_0)$ must hold for any split $S'\in \Sg-\Sg_0$. Thus, putting $\psi:=\phi|_{\Sg_0}$, 
 we have $\psi\in V(\Sigma_0)$ as well as $\phi=\widetilde{\psi}$. So, $\phi \in B(\Sigma_0)$ holds.
 Thus, also
 $\{\phi\in V(\Sigma): \M\cap\Sigma_0 \neq \emptyset\}\subseteq B(\Sigma_0)$ 
and, therefore, also $\{\phi\in V(\Sigma): \M\cap\Sigma_0 \neq \emptyset\}= B(\Sigma_0)$ must hold, as claimed.

\medskip

(v) Next, choosing any fixed map $\phi_0\in B(\Sigma_0)$, we have also
$
B(\Sigma_0)=\{\phi\in V(\Sigma): \Delta(\phi,\phi_0)\subseteq \Sigma_0)\}
$ 
as 
$$
\Delta(\phi,\phi_0)\subseteq \Sigma_0
\iff \phi(S)=\phi_0(S)=A(S\se \Sigma_0)
$$ 
holds for all 
$\phi_0\in B(\Sigma_0)$ and $\phi\in V(\Sigma)$. 
\medskip

(vi) It is also obvious that the map $\phi_{\Sigma_0}$ as defined above
is indeed contained in $B(\Sigma_0)$ 
and that, given 
any map $\phi\in V(\Sg)$, $\phi_{\Sigma_0}$ is indeed the (necessarily unique) map in $B(\Sigma_0)$ with 
$D(\phi, \psi)=D(\phi, \phi_{\Sigma_0}) + 
D( \phi_{\Sigma_0}, \psi)$ 
for all $\psi\in B(\Sigma_0)$ and, hence, the ``gate'' of $\phi$ in $B(\Sigma_0)$. And $\phi_{\Sigma_0}\in cut(\Sg)$ must hold in case  $\phi\notin B(\Sigma_0)$ as it `separates' $\phi$ from $B(\Sigma_0)$. I.e., there must be edges incident with $\phi_{\Sigma_0}$ whose types must be contained in distinct connected components of $\Gamma(\Sg)$, those leading to $\phi$ and those leading to any map $\psi\in B(\Sigma_0)$ distinct from $\phi_{\Sigma_0}$ which must exist as the cardinality of $B(\Sigma_0)$ must be at least $2$ for any connected component $\Sg_0$ of $\Gamma(\Sg)$.

\medskip
(vii) Next, given any two distinct connected components $\Sigma_0,\Sigma_1\in \pi_0(\Sg)$, let $\psi$ denote the map $\psi:=
\phi_{\Sg_0|\Sg_1}: \Sigma\ra \cP(X)$ that maps any split $S\in \Sigma_0$ onto 
$A(S\se \Sigma_1)$ and any other split $S\in \Sg$ onto $A(S\se \Sigma_0)$. It is obvious that the restriction $\psi|_{\Sg_0}: S \mapsto A(S\se \Sigma_1)$ of $\psi$ to $\Sg_0$ is contained in $V(\Sg_0)$ in view of Assertion $(4.6)$. And it is also obvious that $\psi$ coincides with the extension $\widetilde{\psi|_{\Sg_0}}$ of this restriction $\psi|_{\Sg_0}$. So, it is contained in $B(\Sg_0)$. 

It is also the gate $\phi_{\Sg_0}$ in 
$B(\Sg_0)$ of any map $\phi \in B(\Sg_1)$ as, by definition, also $\phi_{\Sg_0}(S) =\phi(S)= A(S\se \Sigma_1)=$ holds for all $S\in\Sg_0$ and 
$\phi_{\Sg_0}(S) =A(S\se \Sigma_0)$ holds for all $S\notin\Sg_0$.

\medskip
(viii) Furthermore,
any map $\phi\in B(\Sg_0)\cap B(\Sg_1)$ must be a cut vertex of $\B(\Sg)$:  
Indeed, (\ref{B00}) implies that there must be edges incident with any such $\phi$ whose types are contained in $\Sg_0$ as well as edges whose types are contained in $\Sg_1$. 
Thus, the assertions $ (B\cap B:1)$ and $ (B\cap B:2)$ are indeed equivalent and follow from $(B\cap B:3)$. Further, if $\phi\in B(\Sg_0)\cap B(\Sg_1)$ holds,
we must have $\phi(S)=A(S\se \Sigma_0)$ for all $S\in \Sg-\Sg_0$ and  $\phi(S)=A(S\se \Sigma_1)$ for all $S\in \Sg-\Sg_1$. So, there can be only one such map, and $A(S\se \Sigma_0)=A(S\se \Sigma_1)$ must hold for all $S\in \Sg-(\Sg_0 \cup \Sg_1)$ in case such a map exists. So, ``$(B\cap B:2)\iff (B\cap B:3)$'' and ``$(B\cap B:2)\Ra (B\cap B:4)$'' holds.

\medskip
(ix) 
And if $(B\cap B:4)$ holds, i.e., if $A(S\se \Sigma_0)=A(S\se \Sigma_1)$ 
holds for all $S\in \Sg-(\Sigma_0\cup\Sigma_1)$, there exists, of course, no split $S=\{A,B\}\in \Sg-(\Sigma_0\cup\Sigma_1)$ with 
$A\cap A_0=\emptyset$ for some set $A_0$ in 
some split $S_0\in \Sigma_0$ and $B\cap B_1=\emptyset$ for some set $B_1$ in some split $S_1\in \Sigma_1$: Indeed, this would imply $A(S\se \Sigma_0)=B\neq A(S\se \Sigma_1)=A$.  

Furthermore, the map $\psi=\phi_{\Sg_0|\Sg_1}$
coincides with the correspondingly defined map $\phi_{\Sg_1|\Sg_0}$. Indeed, Assertion $(B\cap B:4)$ implies that  both of these maps coincide: Both coincide with $A(S\se \Sigma_0)
=A(S\se \Sigma_1)$ on all splits $S$ in $\Sg-(\Sg_0 \cup \Sg_1)$ while they both coincide  
with  $A(S\se \Sigma_1)$ if  $S\in \Sigma_0$ and with $A(S\se \Sigma_9)$  if  $S\in \Sigma_1$. So, $\psi$ must the unique gate in
$B(\Sg_0)$ of every map in $ B(\Sg_1)$ and, simultaneously, the unique gate in
$B(\Sg_1)$ of every map in $ B(\Sg_2)$. Thus,
$B(\Sg_1)\cap B(\Sg_2)\neq \emptyset$ must hold. I.e., $(B\cap B:4)$ implies also $(B\cap B:2)$. So, all of the assertions $(B\cap B:i)\,\, i=1,2,3,4$ must indeed be equivalent, as claimed.

\medskip
(x) Furthermore, as $\B(\Sigma_0)$ is 
$2$-connected (cf.~Corollary \ref{2-connected}), 
it follows from the definition of a block 
that $B(\Sigma_0)$ is a subset of some block $B \in Bl(\Sigma)$. 
However, 
$B(\Sigma_0)$ must, in fact, coincide 
with that block $B$: Indeed, $\phi_0\in B(\Sigma_0)$ and 
$\phi_1\in B$ implies that, by 
the definition of a block, 
$\Gamma_{\phi}(V^{(\phi)},\phi_0)=\Gamma_\phi(V^{(\phi)},\phi_1)$
must hold for all $\phi\in V(\Sg)-\{\phi_0,\phi_1\}$. 
Thus, by Corollary~\ref{2-connected}, 
there must exist some connected component $\Sg'_0$ of the 
incompatibility graph $\Gamma(\Sg)$
of $\Sg$ with 
$\Delta(\phi_0,\phi_1)\subseteq \Sg_0'$.
Moreover, by replacing $\phi_0$ by any 
map $\phi'_0\in B(\Sigma_0)$ with 
$\emptyset\neq \Delta(\phi_0,\phi'_0)$ if necessary, we may also assume that 
$\Delta(\phi_0,\phi_1)\cap \Sg_0\neq\emptyset$ must hold and, therefore, 
$ \Sg_0=\Sg'_0 \supseteq\Delta(\phi_0,\phi_1)$.
This, in turn, implies that also $\phi_1\in B(\Sigma_0)$ 
and, therefore, also $B=B(\Sigma_0)$ must hold.
\medskip

(xi) We clearly have $B(\Sigma_0) \neq B(\Sigma_0')$ for 
any two distinct connected components $\Sigma_0, \Sigma_0'$ 
in $\pi_0(\Sigma)$, 
(as, e.g., $\Delta(\phi,\psi)\subseteq \Sg_0$ 
and $\Delta(\phi',\psi')\subseteq \Sg_0'$ 
must hold for any $\phi,\psi\in B(\Sigma_0)$ 
and $\phi',\psi'\in B(\Sigma'_0)$). Thus,
the map $\Psi$ is injective. 

\medskip

(xii) And finally, 
$\Psi$ is also surjective as there exists, for any 
block $B \in Bl(\Sigma)$, some connected 
component $\Sigma_0 \in \pi_0(\Sigma)$ of $\Gamma(\Sg)$ 
with $B(\Sigma_0)=B$. Indeed, choose any edge 
$e:=\{\phi_1,\phi_2\}$ in $B$, let 
$S$ be the unique split in $\Delta(\phi_1,\phi_2)$, and let 
$\Sigma_0=\Gamma(\Sg,S)$ be the connected component of $\Gamma(\Sg)$ 
that contains the split $S\in\Sg$.
Then, $S\in \Sg^{(\phi_1)}\cap \Sg^{(\phi_2)}$ 
implies, in view of Assertion~(4.5), that 
$\phi_1,\phi_2\in B(\Sg_0)$ must hold. However, 
according to well-known properties of blocks, 
this can be true only if $B=B(\Sg_0)$ holds.
This completes the proof of the theorem.
\epf

\section{The $X$-tree associated with $\Sg$}

Let us recall first that, given any connected simple graph $\Gamma=(V,E)$ 
with vertex set $V$ and edge set $E\subseteq \binom{V}{2}$, 
the graph $\T(\Gamma)$ with vertex set $\V(\Gamma)$ 
the disjoint union $Bl(\Gamma)\amalg V$ of the set 
$Bl(\Gamma)$ of all blocks of $\Gamma$ and the set $V$ of all of its vertices, 
and edge set $\E(\Gamma)$ the set of all pairs $\{B,v\}$ 
with $B \in Bl(\Gamma)$ and $v \in B$ is well-known 
to always be a tree (cf. \cite[Proposition 3.1.2]{D05}).
Clearly, the degree $\deg_{\T(\Gamma)}(v)$ of any 
vertex $v\in V$, considered as a vertex in $\T(\Gamma)$, 
coincides with the number of 
blocks that contain it, and the degree 
$\deg_{\T(\Gamma)}(B)$ of any block $B$ of $\Gamma$, considered as a 
vertex in $\T(\Gamma)$, coincides with the number 
of vertices it contains.

In consequence, continuing with the notation 
introduced so far, every 
collection $\Sigma$ of $X$-splits gives rise to a 
tree $\T(\Sigma):=\T\big(\B(\Sigma)\big)$ 
with vertex set 
$\V(\Sg):=Bl(\Sigma) \,\amalg\, V(\Sigma)$, the disjoint union of 
$Bl(\Sigma)$ and $V(\Sigma)$, and edge 
set the set $\E(\Sg)$ of all pairs $\{B,\phi\}$ 
with $B \in Bl(\Sigma)$ and $\phi \in  B$.
Moreover, 
the canonical labelling map 
$\phi_\Sg: X \ra V(\Sigma): x\mapsto \phi_x$
from $X$ into the $V(\Sigma)$ can also be viewed as a 
labelling map from $X$ into the 
vertex set $\V(\Sigma)$ of  $\T(\Sigma)$. 

Note further that the degree 
$\deg_{\T(\Sg)}(B)$ of any block $B$ of $\B(\Sg)$, considered as a 
vertex in $\T(\Sg)$, coincides with the number 
of vertices it contains. In addition, the degree 
$\deg_{\T(\Sg)}(\phi)$ of any vertex $\phi\in V(\Sigma)$ of $\B(\Sg)$, 
considered as a vertex in $\T(\Sg)$, coincides with the number 
of blocks that contain it and, hence, 
with the cardinality of $\pi_\phi(\M)$, $\pi_\phi(\Sg)$, 
$\pi_\phi(X^{(\phi)})$ as well as of $\pi_\phi(V^{(\phi)})$.
In particular, a vertex $\phi \in V(\Sg)$ 
is of degree larger than $1$ in $\T(\Sg)$ 
if and only if it is a cut vertex of $\B(\Sg)$.

\begin{figure}
\begin{center}
\includegraphics[height=3.5cm, width=8cm]{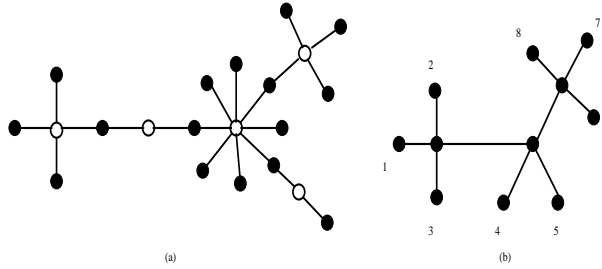}
\caption{(a) The graph $\T(\Sigma_8)$
corresponding to the Buneman graph $\B(\Sigma_8)$ in Figure~\ref{example},
where white vertices correspond to elements in the set $Bl(\Sigma_8)$
and black vertices to those in $V(\Sigma_8)$.
(b)  The $X$-tree $T_{\Sigma_8}$
obtained from the tree $\T(\Sigma_8)$ by
deleting all non-labeled 
vertices $\phi \in \V(\Sg)$ with $\phi \in V(\Sg)-cut(\Sg)$ 
and the pendant 
edges leading to them and
suppressing all vertices of degree $2$ in $\T(\Sg_8)$.
Note that the split system
$\{S_1,S_2,S_3,S_{123},S_4,S_5,S_{678},S_6,S_7,S_8\}$
gives rise to the
same $X$-tree.}\label{construction}
\end{center}
\end{figure}

In consequence, it is fairly obvious that we can also associate, to any 
system $\Sg$ of $X$-splits, a  ``reduced" tree $T_\Sg$
which is a ``proper" $X$-tree: All one needs to do is (i) to 
delete all non-labeled 
vertices $\phi \in \V(\Sg)$ with $\phi \in V(\Sg)-cut(\Sg)$ 
and the pendant 
edges leading to them and (ii) 
to {\em suppress} all vertices of degree $2$, i.e., to replace each maximal 
sequence $u_0,u_1,\dots, u_k$ of distinct non-labeled vertices of 
$\T(\Sg)$ with $\{u_{i-1},u_i\}\in 
\E(\Sg)$ for all $i=1,\dots,k$ and 
$\deg_{T_\Sg}(u_i)=2$ for all $i=1,\dots,k-1$ by just one edge $\{u_{0},u_k\}$ 
while simultaneously deleting all the vertices $u_1,u_2,\dots, u_{k-1}$ 
in between $u_{0}$ and $u_k$ and the 
edges incident with them, see also Figure~\ref{construction}.

Thus, a block $B$ in $Bl(\Sg)$ of the form $B=B(\Sg_0)$ for 
some $\Sigma_0\in \pi_0(\Sg)$ will be suppressed
if and only if  $\Sg_0$ consists of a single split, only.
Otherwise, its degree
$\deg_{T_\Sg}(B)$ in the ``reduced" tree $T_\Sg$ coincides 
with the number 
of equivalence classes of the equivalence 
relation $\sim_{\Sg_0}$ defined on $X$ by 
$$
x \sim_{\Sg_0} y \iff \forall_{S\in \Sg_0} S(x) = S(y).
$$ 
In particular, the labelling map $\phi_\Sg$ sets 
up a bijection from $X$ onto the set of leaves of 
$\T(\Sg)$ if and only if $cut(\Sg)\cap 
\{\phi_x:x\in X\}=\emptyset$  and 
the equivalence relation $\sim_{\Sg}$ is the 
identity relation on $X$, i.e., $\bigcap_{S\in \Sg}S(x)=\{x\}$ 
holds for all $x\in X$.
\medskip

Altogether, this implies 

\begin{theorem}\label{Xtree}
Suppose as above that $X$ is a finite set 
of cardinality at least $2$, and that $\Sigma$ is a 
collection of $X$-splits. Then, the following holds:\\
{\em (i)} The tree $\T(\Sigma)=\T\big(\B(\Sigma)\big)$ canonically associated 
with the graph $\B(\Sigma)$ is isomorphic to the graph with 
vertex set the disjoint union of 
$\pi_0(\Sigma)$ and $V(\Sigma)$, and edge set 
$$
\{ \{\Sigma_0,\phi\} \,:\,   \Sigma_0 \in \pi_0(\Sigma), 
\phi \in V(\Sigma),  \M \cap \Sigma_0 \neq \emptyset \}.
$$
{\em (ii)}  The graph 
obtained from the tree $\T(\Sigma)$ by deleting  all non-labeled 
vertices $\phi \in  V(\Sg)-cut(\Sg)$ 
and the pendant 
edges leading to them as well as
``suppressing'' all non-labeled vertices of degree $2$ in $\T(\Sg)$  is, 
together with the induced labelling map from $X$ 
into its vertex set, an $X$-tree that we denote by $T_\Sg$. 
It is canonically associated with $\Sg$ and coincides 
with that $X$-tree that, according to Peter Buneman, 
is associated with $\Sg$ in case $\Sg$ is compatible.

{\em (iii)}
Furthermore, given any distinct 
blocks $B_1, B_2\in Bl(\Sigma)$ for which a 
cut vertex $\phi\in B_1 \cap B_2$ exists, one has 
$$
\max\big(\deg_{\T(\Sg)}(B_1),\deg_{\T(\Sg)}(B_2),\deg_{\T(\Sg)}(\phi) 
\big) \ge 3,
$$ 
that is, at least one of the three vertices 
$B_1,B_2,$ and $\phi$ of $\T(\Sg)$ 
must, for every such triple $B_1,B_2,$ 
and $\phi$, also be a vertex of 
the $X$-tree $T_\Sg$ derived from $\T(\Sg)$.
\end{theorem}
\pf
The first two assertions follow immediately from 
our previous observations. To establish (iii), note
that since $B_1\not= B_2$ we have 
$\deg_{\T(\Sg)}(B_1),\deg_{\T(\Sg)}(B_2)\ge 2$. But 
$\deg_{\T(\Sg)}(B_1),\deg_{\T(\Sg)}(B_2)= 2$ would imply that 
there exist two distinct vertices $\psi_1,\psi_2 \in V(\Sigma)-\{\phi\}$
with $B_1=\{\psi_1,\phi\}$ and 
$B_2=\{\psi_2,\phi\}$.
Let $S_1,S_2\in \Sg$ denote the splits for which 
$\Delta(\psi_1,\phi)=\{S_1\}$ and $\Delta(\psi_2,\phi)=\{S_2\}$ 
both hold. Then $A_1 \cup A_2 =X$ 
must hold for $A_1:=\phi(S_1)$  
and $A_2:=\phi(S_2)$. This, in turn, implies that 
$A_1\cap A_2 \neq \emptyset$ must hold as well as
$S_1,S_2\not\in\Delta(\phi,\phi_x)$
for every $x\in A_1\cap A_2$. Hence,
in view of Lemma~\ref{incomp} and Theorem~\ref{blockbijection},
$\Gamma_\psi(\M)$ contains at least three connected
components: $\{S_1\}, \{S_2\}$, and the connected
component containing $\Delta_{\min}(\phi|\phi_x)$
for some $x\in A_1\cap A_2$.
\epf

As mentioned already above, we will 
explore these matters in more detail in \cite{DHKM07e}. We will consider 
in particular
the graph theoretical invariant, defined according to \cite{DHKM07c}, of 
the {\em Buneman complex} 
${\bf B}(\Sigma)$ associated to $\Sigma$ 
as defined in \cite{dre-hub-mou-97}.
And we will show that $\B(\Sg)$ as well as  ${\bf B}(\Sigma)$
can be described as a push-out
in terms of the data involved in the construction of  $T_\Sg$ suggesting 
efficient algorithms for their computation (see also \cite{GB05, GB07}).

\medskip
In yet another paper, we will discuss what can be done in case the 
incompatibility graph $\Gamma(\Sg)$ is connected and $\B(\Sigma)$ is, 
hence, 2-connected. In particular, by directing attention 
towards {\em cut faces} 
rather than merely cut vertices of $\B(\Sigma)$, 
we show that it may be possible 
even in this case to extract valuable phylogenetic 
information from a split system $\Sg$.\\

\noindent{\bf Acknowledgment}:
The authors thank 
the Engineering and Physical Sciences
Research Council (EPSRC) for its support [Grant EP/D068800/1], and also 
Andreas Spillner for helpful discussions.
K.\,T.\,Huber and J.\,H.\,Koolen thank the Royal Society for its support
through their International Joint Projects scheme. 
J.\,H.\,Koolen was also partially supported by the Priority 
Research Centers Program
through the National Research Foundation of 
Korea (NRF) funded by the Ministry of
Education, Science and Technology (Grant 2009-0094069).
A.\,Dress thanks the Chinese Academy for Sciences,
the Max-Planck-Gesellschaft, and the German BMBF for their support, as
well as the Warwick Institute for Advanced Study where,
during two wonderful weeks, the basic outline of this paper was conceived.
Finally, we also thank an anonymous referee for his/her 
helpful comments on the first version of this paper.

\end{document}